\documentclass[AMA,STIX1COL]{WileyNJD-v2}

\articletype{Research Article}%

\received{26 April 2016}
\revised{6 June 2016}
\accepted{6 June 2016}

\newcommand{\bunderline}[1]{\underline{#1\mkern-4mu}\mkern4mu }
\newcommand{\overbar}[1]{\mkern 1.5mu\overline{\mkern-1.5mu#1\mkern-1.5mu}\mkern 1.5mu}
\newcommand{\fstar}{f^\ast}
\newcommand{\knom}{k_1} %{k_\mathrm{nom}}
\newcommand{\kBackstep}{k_2} %{c}
\newcommand{\kSafe}{k_3} %{k}
\newcommand{\Dstar}{D^\ast }
\newcommand{\ystar}{y^\ast }

\newcommand{\Dnom}{D_\mathrm{nom}}
\newcommand{\LyapCoeff}[1]{b_{#1}}
\newcommand{\KarafControlGain}[1]{k_{#1}}
\newcommand{\uCancel}{u_\mathrm{c}}
\newcommand{\uStabilize}{u_\mathrm{s}}
\newcommand{\realNumbers}{\mathbb R}
\newcommand{\Dmin}{\bunderline D}
\newcommand{\Dmax}{\overbar D}
\newcommand{\KarafParA}{a_1} % Miroslavs suggestion:\Delta D

\newcommand{\KarafParB}{a_2} % Miroslavs suggestion {\KarafParBmax}/{\KarafParBmin}

\newcommand{\diff}{{\rm d}}

\usepackage{ifthen} % for \exp
\renewcommand{\exp}[1][]{%
   \ifthenelse{ \equal{#1}{} }
      {{\mathrm{exp}}}
      {{\mathrm{exp}\left( #1 \right)}}
}

\usepackage{mathtools}

\makeatletter
\newcommand*{\citenumns}[2][]{%
  \begingroup
  \let\NAT@mbox=\mbox
  \let\@cite\NAT@citenum
  \let\NAT@space\NAT@spacechar
  \let\NAT@super@kern\relax
  \renewcommand\NAT@open{}%
  \renewcommand\NAT@close{}%
  \cite[#1]{#2}%
  \endgroup
}
\makeatother
\begin{document}

\title{Stabilization of Age-Structured Chemostat Hyperbolic PDE with Actuator Dynamics}
%and Delay}
% \title{This is the sample article title\thanks{This is an example for title footnote.}}

\author[1]{Paul-Erik Haacker}

\author[2]{Iasson Karafyllis}

\author[3]{Miroslav Krstić}

\author[3]{Mamadou Diagne*}

\authormark{Haacker \textsc{et al}}

\address[1]{\orgdiv{Institute for Nonlinear Mechanics}, \orgname{University of Stuttgart}, \orgaddress{\state{Pfaffenwaldring 9, 70569 Stuttgart}, \country{Germany}}}

\address[2]{\orgdiv{Department of Mathematics}, \orgname{National Technical University of Athens}, \orgaddress{\state{Zografou Campus, 15780 Athens}, \country{Greece}}}

\address[3]{\orgdiv{Department of Mechanical and Aerospace Engineering}, \orgname{University of California San Diego}, \orgaddress{\state{La Jolla, CA 92093-0411}, \country{USA}}}

\corres{*Mamadou Diagne, \orgdiv{Department of Mechanical and Aerospace Engineering}, \orgname{University of California San Diego}, \orgaddress{\state{La Jolla, CA 92093-0411}, \country{USA}}. \email{mdiagne@ucsd.edu}}

% \presentaddress{This is sample for present address text this is sample for present address text}

\abstract[Abstract]{
%\color{magenta}
For population systems modeled by age-structured hyperbolic partial differential equations (PDEs), we redesign the existing feedback laws, designed under the assumption that the dilution input is directly actuated, to the more realistic case where dilution is governed by actuation dynamics (modeled simply by an integrator). In addition to the standard constraint that the population density must remain positive, the dilution dynamics introduce constraints of not only positivity of dilution, but possibly of given positive lower and upper bounds on dilution. We present several designs, of varying complexity, and with various measurement requirements, which not only ensure global asymptotic (and local exponential) stabilization of a desired positive population density profile from all positive initial conditions, but do so without violating the constraints on the dilution state. 
To develop the results, we exploit the relation between first-order hyperbolic PDEs and an equivalent representation in which a scalar input-driven mode is decoupled from input-free infinite-dimensional internal dynamics represented by an integral delay system.}


%\jnlcitation{\cname{%
%\author{Haacker P.-E.},
%\author{I. Karafyllis},
%\author{M. Krstić}, and
%\author{M. Diagne}} (\cyear{2023}),
%\ctitle{Stabilization of Age-Structured Chemostat Hyperbolic PDE with Actuator Dynamics}, \cjournal{Int J Robust Nonlinear Control}, \cvol{???}.}

\maketitle

\footnotetext{\textbf{Abbreviations:} ODE, Ordinary Differential Equation; PDE, Partial Differential Equation; IDE, Integral Delay Equation; BC, Boundary Condition; OE, Output Equation}

\section{Introduction}

\subsection*{Motivation: Chemostat as a benchmark for nonlinear control and for control in epidemiology} 

Many industrial biotechnology processes can be described by a nonlinear PDE (Partial Differential Equation) of population dynamics that are  structured around age. Reflecting the maturation of enzymes, microbial organisms, animal or plant cells and tissues \cite{kuritz2017relationship,billy2012age,rong2007mathematical}, these models of biological and biochemical systems are often exploited  to favor the desired functionalities  and achieve cost-effective  bioreactors' production rate \cite{KK17Chemostat,  kurth2021tracking,kurth2023control}. The word ``chemostat'' denotes a biological process which is fed a nutrient at a certain rate and from which a bioproduct-nutrient mix is being extracted/removed at the same rate \cite{KK17Chemostat}. This rate is referred to as the ``dilution rate.'' The relevance of chemostats in control engineering extends beyond their own application in biotechnology (such as the manufacturing of  pharmaceuticals). First, a chemostat is a nonlinear control problem---even when the limiting substrate(s) nonlinear dynamics are neglected---due to the fact that its dynamic model contains a product of the population density and the dilution rate. Second, the relevance of the chemostat is also in the fact that it is an example of a ``positive system,'' namely, a control system whose state is subject to an inequality constraint \cite{de2001stabilization2,de2001stabilization}. 

Third, and perhaps most important at present, chemostat is a special case of more complex population dynamics such as those that arise in epidemiology \cite{inaba2017age,martcheva2015introduction}. Before one takes on control design of SIR-type (susceptible-infected-recovered) models in epidemiology, in which control is being applied to more than one population (for example, to the susceptible and the infected populations), it is necessary to understand how to control single-population chemostats. Just as there is a parallel between the biomass in a chemostat and the infected population in epidemiology, as well as between the substrate in a chemostat and the susceptible population in epidemiology \cite{toth2006limit,inaba2017age}, there is likewise a parallel between the treatment/therapy being administered as an input in epidemiology and the dilution input of a chemostat. 

Hence, the study of control of an age-structured chemostat that we undertake here is an important first step towards designing controllers for {\emph{age-structured}} models in epidemiology.

\subsection*{Age-structured population dynamics} 

Age-structured chemostat models are described by a particular first order hyperbolic Partial Differential Equation (PDE) with a non-local boundary condition: the McKendrick-von Foerster equation (see References \citenumns{m1925applications,brauer2012mathematical,kermack1927contribution,boucekkine2011optimalControl,charlesworth1994ageStructured,rundnicki1994asymptoticSimilarity} and the references therein). Age-structured models are extensions of the chemostat models described by Ordinary Differential Equations (ODEs; see Reference \citenumns{smith1995theoryOfChemostat}). A particular and deep mathematical tool has been developed for the study of age-structured models: the ergodic theorem, also known as the asymptotic behavior (see References \citenumns{inaba1988semigroupErgodic,inaba1988asymptoticProperties, rundnicki1994asymptoticSimilarity} for similar results on asymptotic similarity and  Reference \citenumns{KK17Chemostat} for the proof of the ergodic theorem by means of Lyapunov arguments). 

Control studies for age-structured models are rare. Optimal control problems have been studied (see References \citenumns{boucekkine2011optimalControl,feichtinger2003optimality,sun2014optimal} and the references therein). On the other hand, feedback control of infinite-dimensional population dynamics was introduced in Reference \citenumns{KK17Chemostat} and further considered in References \citenumns{schmidt2018trajectory,kurth2023control}. More specifically, the use of the ergodic theorem and the corresponding study of linear integral delay equations in References \citenumns{KK17Chemostat} led to a special nonlinear infinite-dimensional change of variables, from the age-structured population density (a population density that is a function of the continuum age) into 
\begin{itemize}
    \item[(i)] a scalar variable that represents the controllable mode of the bilinear non-local hyperbolic PDE system which is directly actuated by the dilution rate input, and
\item[(ii)] an infinite-dimensional uncontrollable but exponentially stable dynamical system described by linear integral delay equations,
\end{itemize}
 and has subsequently been the cornerstone of the control design and stability analysis of age-structured population dynamics. 

The control studies in References \citenumns{KK17Chemostat} and \citenumns{schmidt2018trajectory,kurth2023control} used the dilution rate as the control input as in many other control studies of chemostat models described by ODEs (see References \citenumns{karafyllis2008controlODE,karafyllis2009relaxedODEcontrol,gouze2006robust,de2003feedback,dimitrova2012nonlinear}). Actuation by dilution rate is both physical and plausible and it amounts to harvesting the population. In simplistic terms, a steady level of harvesting is needed to maintain the population at a sustained level (equilibrium) and, when the population exceeds such an equilibrium, it is over-harvested by (over)dilution, whereas when the population drops below such an equilibrium, it is under-harvested by (under)dilution. How this over- and under-harvesting is to be done exactly is highly non-trivial because the population is age-structured, meaning that its state is functional/infinite-dimensional, due to which one cannot really speak of overpopulation and underpopulation in a bulk sense. One has, instead, to take into account that the population can vary from the equilibrium in infinitely many ways (overpopulated old with underpopulated young, vice-versa, and so on). This infinite variety in the state profiles, in addition to the infinite number of states (age-specific densities), along with the nonlinearity of the problem, is one of the theoretical attractions to this control design problem.

\subsection*{Chemostat with dynamic actuation of dilution} 

While dilution is a plausible form of actuation, it is not possible to actuate it instantaneously. As with any other actuation, it is governed by its own dynamics. In the case of dilution, it is affected by the dynamics of pipes and valves that are used to run the nutrients into a chemostat and the population-nutrient mix out of the chemostat. In biochemical engineering, the dilution rate is defined as the ratio of the inlet volumetric flow rate to the reactor volume. Since the volume of the reactor is not necessarily constant, it follows that the dynamics of the dilution rate can play an important role. 

In this paper we design stabilizing feedback laws for the age-structured chemostat with dilution dynamics governed by a single, scalar integrator. 
%two different simple models: (i) an integrator, \textcolor{red}{ (ii) a pure delay.} 
In the most advanced versions of our designs, we ensure that not only the population density, but also the dilution rate, remains in a specific interval (and positive). Dilution rate being positive means that we are not feeding population back into the chemostat but are only harvesting/removing population. 

While the age-structured chemostat with a constant dilution rate is a PDE system with one eigenvalue at the origin and all other eigenvalues in the left half plane, meaning that it is not open-loop unstable, and meaning that it can, at worst, settle at an undesirable age profile for the population (provided the constant dilution is at a population-sustaining level), the system with dilution dynamics modeled by an integrator has two eigenvalues at the origin, meaning that it is open-loop unstable. This is one more way of noting the need for a more sophisticated control design for this PDE system once the dynamics of the dilution actuation are introduced. 

In this paper we tackle three design challenges:

\begin{itemize}
    \item[(i)] we overcome the dynamics of the integrator modeling the dilution actuation,
    \item[(ii)] we perform the design and the analysis using feedback of only a sensor of bulk population, which has an unknown profile in sensing populations of different ages,
    \item[(iii)] we ensure that the dilution rate that drives the population to the desired age-structured equilibrium from all initial positive population density age profiles remains itself in a specific interval (and positive).
%   \item[(iv)] \textcolor{red} {we solve the feedback stabilization problem when the dynamics of the dilution actuation are described by a single delay.}
\end{itemize}

While in the initial work on control of the age-structured chemostat \cite{KK17Chemostat} the state space of the systems was, so to speak, the ``positive orthant'' of infinite dimension (the age-structured population density), in this work the state space increases by one positive-valued state when the dilution dynamics are modeled by an integrator and becomes, so to speak, a ``positive orthant'' of dimension $\infty+1$. From a mathematical point of view, the paper deals with a global feedback stabilization problem for a constrained nonlinear infinite-dimensional control system which is defined on a specific open set. Finite-dimensional control systems defined on open sets were recently studied in Reference \citenumns{sontag2022remarks}. 

\bigskip

\textbf{\emph{Organization: }}
The structure of the paper is as follows. Section \ref{sec_popDynamicsActuatorDynamicsBackstepping-model} introduces the mathematical model as well as the stabilization problem. Sections \ref{sec_popDynamicsActuatorDynamicsBackstepping} and \ref{sec_BacksteppingStabilityProof} study the backstepping design of the problem with unrestricted dilution rate assuming that the state is available, while Section \ref{sec_RelaxedBackstepping} provides globally stabilizing output feedback laws. Sections \ref{sec_Safety} and \ref{sec-constrained-miroslav} are devoted to the feedback stabilization problem of the more demanding case where the dilution rate is restricted to take values in a specific interval. 
%Finally, Section 8 studies the design of predictor feedback controllers in the case where there is a delay in the dilution rate. 

\textbf{\emph{Notation:}}

%\noindent \textbf{Notation.} Throughout this paper, we adopt the following notation. 

\begin{itemize}

\item[*]  $\realNumbers_+ $ denotes the interval $[0,+\infty )$.

\item[*] Let $U$ be an open subset of a metric space and $\Omega \subseteq \realNumbers^{m} $ be a set. By $C^{0} (U ; \Omega )$, we denote the class of continuous mappings on $U$, which take values in $\Omega $. When $U\subseteq \realNumbers^{n} $, by $C^{1} (U ; \Omega )$, we denote the class of continuously differentiable functions on $U$, which take values in $\Omega $. When $U=[a,b)\subseteq \realNumbers $ (or $U=[a,b]\subseteq \realNumbers $) with $a<b$, $C^{0} ([a,b) ; \Omega )$ (or $C^{0} ([a,b] ; \Omega )$) denotes all functions $f:[a,b)\to \Omega $ (or $f:[a,b]\to \Omega $), which are continuous on $(a,b)$ and satisfy $\mathop{\lim }\limits_{s\to a^{+} } \left(f(s)\right)=f(a)$ (or $\mathop{\lim }\limits_{s\to a^{+} } \left(f(s)\right)=f(a)$ and $\mathop{\lim }\limits_{s\to b^{-} } \left(f(s)\right)=f(b)$). When $U=[a,b)\subseteq \realNumbers $, $C^{1} ([a,b) ; \Omega )$ denotes all functions $f:[a,b)\to \Omega $ which are continuously differentiable on $(a,b)$ and satisfy $\mathop{\lim }\limits_{s\to a^{+} } \left(f(s)\right)=f(a)$ and $\mathop{\lim }\limits_{h\to 0^{+} } h^{-1} \left(f(a+h)-f(a)\right)=\mathop{\lim }\limits_{s\to a^{+} } f'(s)$. 

\item[*] $\mathcal K_{\infty} $ is the class of all strictly increasing, unbounded functions $a\in C^{0} ( \realNumbers_+ ; \realNumbers_+ )$, with $a(0)=0$ (see Reference \citenumns{karafyllisJiang2011stability}).

\item[*] $\mathcal{KL}$ is the class of functions $\beta:\realNumbers_+ \times \realNumbers_+ \to \realNumbers_+$ which satisfy the following: For each $t\geq 0$, the mapping $\beta(\cdot,t)$ is of class $\mathcal{K}$, and, for each $s\geq 0$, the mapping $\beta(s,\cdot)$ is nonincreasing with $\lim_{t \to \infty} \beta(s,t) = 0$ (see Reference \citenumns{karafyllisJiang2011stability}).

\item[*] For any subset $S\subseteq \realNumbers $ and for any $A>0$, $PC^{1} \left([0,A];S\right)$ denotes the class of all functions $f\in C^{0} ([0,A];S)$ for which there exists a finite (or empty) set $B\subset (0,A)$ such that: (i) the derivative $f'(a)$ exists at every $a\in (0,A)\backslash B$ and is a continuous function on $(0,A)\backslash B$, (ii) all meaningful right and left limits of $f'(a)$ when $a$ tends to a point in $B\bigcup \{ 0,A\} $ exist and are finite. 

\item[*] Let a function $f\in C^{0} (\realNumbers_+ \times [0,A])$ be given, where $A>0$ is a constant. We use the notation $f[t]$ to denote the profile at certain $t\ge 0$, i.e., $(f[t])(a)=f(a,t)$ for all $a\in [0,A]$.

\item[*] Let a function $x\in C^{0} ([-A,+\infty );\realNumbers )$ be given, where $A>0$ is a constant. We use the notation $x_{t} \in C^{0} ([-A,0];\realNumbers )$ to denote the ``$A-$history'' of $x$ at certain $t\ge 0$, i.e., $\left(x_{t} \right)(-a)=x(t-a)$ for all $a\in [0,A]$.

\item[*] Let a function $f\in C^0(\realNumbers_{+}; \realNumbers)$ be given. By $D^+ f:\realNumbers_+ \to \realNumbers$, \newline$t \mapsto (D^+ f)(t) = \limsup_{h\to 0^+}\left\{h^{-1}(f(t+h)-f(t))\right\}$, we denote the upper Dini derivative.

% \item[*] The inner product of $L^2$ is denoted $\langle f,g\rangle \coloneqq \int_0^A f(a)g(a)\diff a$ where $f,g\in L^2([0,A])$.

\item[*] Let $A,\;B$ be logical statements. Their logical conjunction is denoted by $A\land B$.

% upper right Dini-derivative
%     \begin{equation}
%         D^+ f(t) = \limsup_{h\to 0^+}\left\{h^{-1}(f(t+h)-f(t))\right\}
%     \end{equation}
%     with 

\end{itemize}

\section{Population Model with Actuator Dynamics}
\label{sec_popDynamicsActuatorDynamicsBackstepping-model}

We consider   the age distribution of the population of a microorganism in a bioreactor governed by   the following equations 
\begin{subequations}
    \label{eq_systemOpenLoop}
    \begin{alignat}{3}
         \text{PDE:}& \quad &&\frac{\partial {f}}{\partial t}(a,t) +\frac{\partial f }{\partial a}(a,t)=-\left(\mu(a)+D(t)\right)f(a,t) \label{eq_PDE}\\
         \text{ODE:}& \quad &&\dot D (t)=u(t) \label{eq_ODE}\\
         \text{OE:}& \quad &&y(t) = \int_0^A p(a)f(a,t) \diff a\label{eq_output}\\
         \text{BC:}& \quad &&f(0,t)=\int_0^A k(a)f(a,t) \diff a \label{eq_OLBC}%\\
 %        \text{IC:}& \quad &&f(a,0) = f_0(a)\geq 0, \quad a \in [0.A]\label{eq_IC} 
    \end{alignat}
\end{subequations}
{where $a \in [0,A]$ is the age variable, $A$ is the maximum reproductive age  of the microorganism (a constant), $t\in \realNumbers_+$ is time, $f(a,t)>0$ denotes the distribution of the microbial mass in the reactor at age $a \in [0,A]$ and time $t$, $\mu$ and $k$ are the mortality rate profile and birth rate profile, respectively, 
%continuous in age $a$, 
and $p$ is the sensor kernel. The functions $\mu :[0,A]\to \realNumbers_+ $ and $k:[0,A]\to \realNumbers_+ $ are assumed to be continuous in age $a$, and the function $p :[0,A]\to \realNumbers_+ $ continuously differentiable in $a$ with $\int_0^Ap(a)\diff a >0$.
%, continuously differentiable in age $a$. 
The population of microorganisms grows at a rate regulated  by  the dilution rate  $D(t),$ which is  the control input of model \eqref{eq_systemOpenLoop}. Accounting for  the inlet and outlet rate of the bioreactor, the PDE \eqref{eq_PDE} expresses the microbial mass balance as a variation of the \textit{McKendrick-von Foerster} equation \cite{m1925applications,kermack1927contribution,brauer2012mathematical,charlesworth1994ageStructured,rundnicki1994asymptoticSimilarity}. The internal boundary feedback \eqref{eq_OLBC}, is the renewal condition, which expresses the reproduction of the microorganism as  the current mass of the newborns $f(0,t)$. From a mathematical point of view, the boundary condition \eqref{eq_OLBC} involves  non-local terms. The measured output defined in \eqref{eq_output}  is a weighted average of the mass of the microorganism with an age-specific kernel $p(a)$. }

{ In chemostat reactors,  the dilution rate is defined as  the ratio of the inlet volumetric flow rate to the reactor volume. Although,    the reactor's inlet flow  and outlet volumetric flow rates are  readily adjustable, the volume of the  reactor is not necessarily constant and its variations can be described by   the dynamic actuation introduced as \eqref{eq_ODE}. Indeed,  the population balance \eqref{eq_systemOpenLoop} is a realistic extension of the model studied in Reference \citenumns{KK17Chemostat}, which does not account for a  dynamic actuation and consequently  assumes no possible variations of   the reactor volume. Equation \eqref{eq_ODE} defines an appropriate control input $u(t) \in \realNumbers$ that depends on the reactor volume, the reactor inlet and outlet volumetric flow rates as well as on the time derivative of the reactor inlet volumetric flow rate. }

{It should also be noticed that model \eqref{eq_systemOpenLoop} is  derived by neglecting the dependence of the growth of the microorganism on the concentration of a limiting substrate. However, resource-based models that incorporate limiting substrates or uptake of  nutrients are more suitable to describe the dynamics of continuous microbial culture and might exhibit limit cycles that are  induced by the behavior of a limiting resource assuming constant dilution rates \cite{toth2006limit,toth2008strong}.  Model \eqref{eq_systemOpenLoop} presupposes  that the growth or decay of  a living population only depends on  irreversible incidences of birth and death and is therefore suitable to predict the evolution of  macro-populations in  demography, epidemiology \cite{inaba2017age} and ecology \cite{toth2008bifurcation}. }

%{\color{red}and $p \in C^1([0,A])$}.  
The state of model \eqref{eq_systemOpenLoop} is $(f[t],D(t)) \in \mathcal{F} \times \realNumbers$, where $\mathcal{F}$ is the function space defined by the following equation: 
\begin{align}
    \mathcal{F} = \left\{ f \in PC^1([0,A];(0,\infty)): f(0)=\int_0^A k(a)f(a) \diff a\right\}.
\end{align}
{The results of this contribution encompass   two different cases for the state space. For the first  case, the state space is defined as  $\mathcal{F} \times \realNumbers$, allowing for the dilution rate $D(t)$  to take arbitrary real values. However, as noted above, the dilution rate is the ratio of the inlet volumetric flow rate to the reactor volume and consequently, it is physically meaningful to restrict its set of admissible values to the interval $I \subseteq \realNumbers_+$, which is  much more demanding as the state is constrained and the state space defined as $\mathcal{F} \times I$. }

Our main assumption for model \eqref{eq_systemOpenLoop} is the existence of a constant $\Dstar>0$ such that 
\begin{align} 
1=\int _{0}^{A}k(a)\exp \left(-D^{*} a-\int _{0}^{a}\mu (s)ds \right)da,  
\end{align} 
{which is the analogue of the Lotka-Sharpe equation (see Reference \citenumns{brauer2012mathematical})  and equivalent to assuming population viability. When the above assumption holds, model \eqref{eq_systemOpenLoop} admits a family of equilibrium points $(\fstar,\Dstar)$ (as the model without actuator dynamics \cite{KK17Chemostat}), given by the following equation for arbitrary $M>0$:}
\begin{align}
\fstar(a) = M \exp \left(-D^{*} a-\int _{0}^{a}\mu (s)ds \right). 
\end{align}
\label{eq_steadyState}
{The existence of a continuum of equilibria for model \eqref{eq_systemOpenLoop} implies that none of the equilibrium points is asymptotically stable. Therefore, if we want the bioreactor to operate on a specific equilibrium point, then we need to stabilize such an equilibrium point by means of feedback control. }

\section{Backstepping Design with Unrestricted Dilution}
\label{sec_popDynamicsActuatorDynamicsBackstepping}
%\subsection{Backstepping Control Design}

{Our objective is  to stabilize a desired steady-state \eqref{eq_steadyState} for a chosen $M>0$ or equivalently, stabilize an appropriate output setpoint 
\begin{align}\label{eq_ystar_def}
    \ystar = \int_0^A p(a) \fstar(a) \ \diff a >0.
\end{align} Our design is build upon the control law proposed  in Reference \citenumns{KK17Chemostat} for the stabilization   of system \eqref{eq_systemOpenLoop} without actuator dynamics. We recall  the output-feedback controller  \cite{KK17Chemostat}
\begin{equation}
    \Dnom (t) = \Dstar + \knom\ln \frac{y(t)}{\ystar},~~\knom >0
    \label{eq_nominalControl}
\end{equation} that stabilizes the equilibrium corresponding to $\ystar$. To extend controller \eqref{eq_nominalControl} to the case with actuator dynamics described by system \eqref{eq_systemOpenLoop},   a backstepping approach is adopted. Defining the dilution rate error as follows
\begin{align}
    \delta(t) \coloneqq D(t) - \Dnom (t) \label{eq_defintionDelta}
\end{align}
and taking its time-derivative along the output $y(t)$ given by \eqref{eq_output}, with the help of integration by parts, the following ODE is derived
\begin{align}
    \dot \delta (t) &= u(t) - \knom \frac{\dot y(t)}{y(t)}\\
    &= u(t) + \knom D(t) + \frac{\knom}{y(t)}\left[ p(A)f(A,t) - p(0) f(0,t) - \int_0^A \Tilde{p}(a) f(a,t) \diff a \right]
\end{align}
where 
%we assume $p \in C^1([0,A])$ in Sections \ref{sec_popDynamicsActuatorDynamicsBackstepping},\ref{sec_BacksteppingStabilityProof} and
% \ref{sec_RelaxedBackstepping}
\begin{equation}
    \Tilde{p}(a) \coloneqq p'(a) -p(a)\mu(a).
\end{equation}
At this point, we aim at designing a control law $u(t)$ that achieves exponential convergence of the error $\delta (t)$. Selecting $k_2>0$ such that
\begin{equation}
    \dot\delta (t) = -\kBackstep \delta(t),
\end{equation} the following  controller is designed
\begin{subequations}
    \label{eq_controllerBackstepping}
    \begin{align}
        u(t) &= \uCancel(t) + \uStabilize(t)\\
        \uCancel(t) &= - \knom D(t) - \frac{\knom}{y(t)}\left[ p(A)f(A,t) - p(0) f(0,t) - \int_0^A \Tilde{p}(a) f(a,t) \diff a \right]\label{eq_cancelinggeneral}\\
        \uStabilize(t) &= -\kBackstep \left( D(t) - \Dstar -\knom \ln \frac{y(t)}{\ystar}\right) = -\kBackstep \delta(t),
    \end{align}
\end{subequations}
 which is a  superposition of canceling and stabilizing terms $\uCancel(t)$ and $\uStabilize(t)$. In the most general case, the controller \eqref{eq_controllerBackstepping} requires full measurement of the plant and the actuator states $f(a,t)$ and $D(t),$ respectively. Furthermore, the feasibility of controller \eqref{eq_controllerBackstepping} depends on   the knowledge of the parameters  $p$, $p'$, $\mu$ and $\Dstar$. However,  assuming constant kernel $p$ and mortality rate $\mu$, the  canceling terms \eqref{eq_cancelinggeneral} of backstepping controller \eqref{eq_controllerBackstepping} reduces to the following feedback law   
\begin{equation}
    \uCancel(t) = - D(t) - \frac{p}{y(t)}\left[ f(A,t) - f(0,t) \right] - \mu,
\end{equation}
which removes the burden of full-state measurement of the plant's state.}

\begin{figure}[t]
    \centering
    \includegraphics[width=\textwidth]{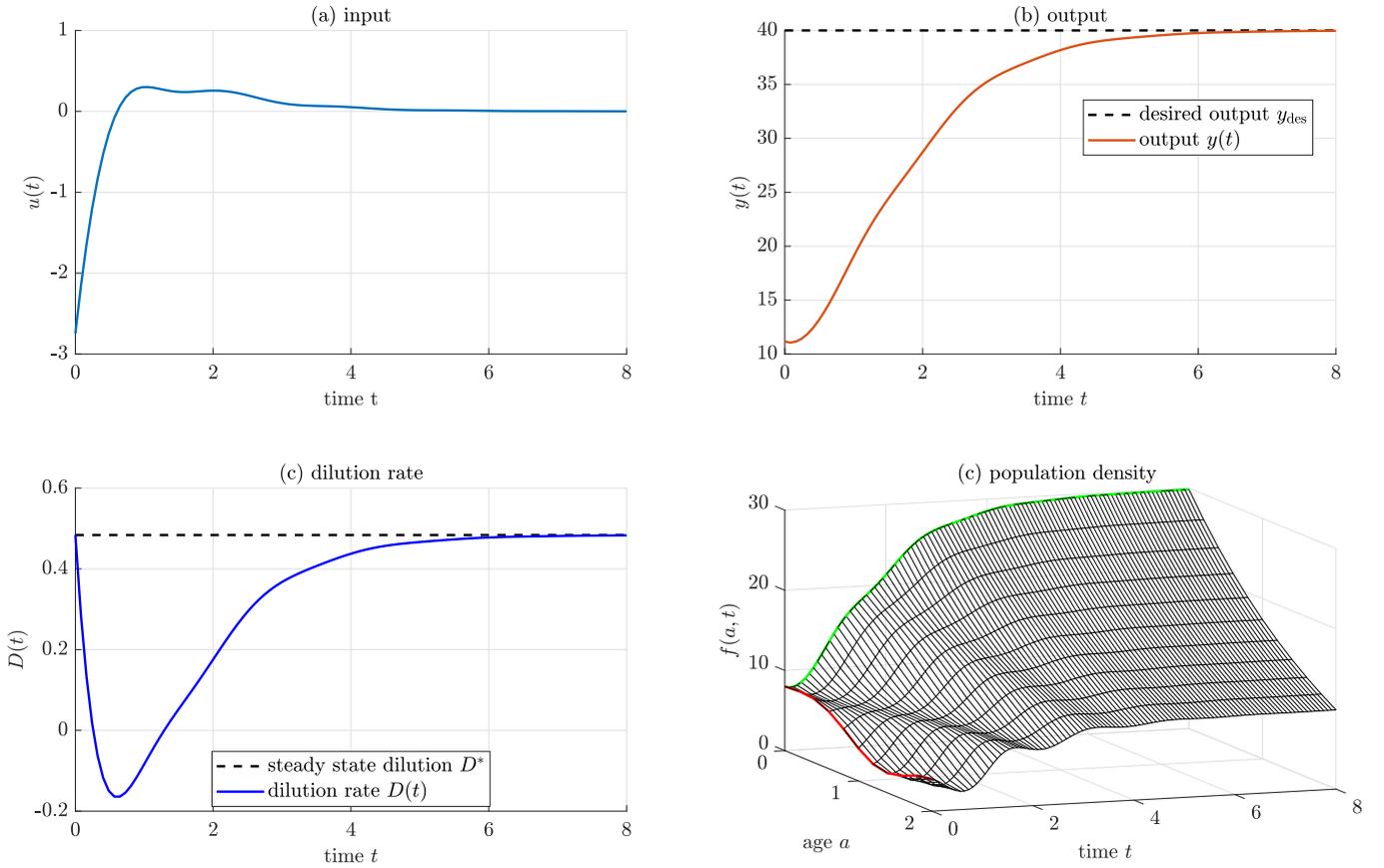}
    % \caption{Simulation results of system \eqref{eq_systemOpenLoop} with controller \eqref{eq_controllerBackstepping} and parameters $\mu = 0.2$, $p = 1$, $k(a) = 2a(A-a)$ and $A = 2$. The controller acts without restricting the dilution rate $D(t)$ to physically meaningful values.}
    \caption{Simulation results of system \eqref{eq_systemOpenLoop} under the control law \eqref{eq_controllerBackstepping} and with the parameters and the initial condition given by \eqref{eq_SimulationParameterAndInitialCondition}. The  gains of the controller are chosen as $\KarafControlGain{1} = 1$ and $\KarafControlGain{2} = 2.$  Convergence to the steady-state is achieved but with possibly negative dilution rate $D(t).$}
    \label{fig_ActuatorDynamics_staticController_AllFinal}
\end{figure}

{
To illustrate the stabilizing effect of the control law \eqref{eq_controllerBackstepping}, we first present simulation results before conducting a stability analysis. Following Galerkin's approach \cite{kurth2021optimalControl,schmidt2017Simulation},  the approximate population density, $\hat f: [0,A]\times \realNumbers_+ \to \realNumbers$ is defined as follows
\begin{align}
    \hat f(a,t) = \sum_{j=1}^N \phi_j(a)\xi_j(t) \eqqcolon \phi(a)^\top \xi(t)
\end{align}
where  $\phi_j: [0,A]\to \realNumbers$ are the trial functions satisfying the boundary condition \eqref{eq_OLBC}, chosen as in Reference \citenumns{kurth2021optimalControl} with $N=6$, and the temporal weights $\xi_j:\realNumbers_+ \to \realNumbers$  are solutions to the initial value problem below
\begin{subequations}
    \begin{align}
    \dot \xi(t) &= (A_\phi - D(t) I)\xi(t), \quad t\geq 0\\
    \xi(0) &= \xi_0
    \end{align}
\end{subequations}
where $A_\phi \in \realNumbers^{N\times N}$ is determined by the choice of trial functions $\phi_j$, and $\xi_0 \in \realNumbers^{N}$ is computed for a given initial population density $f_0 \in \mathcal{F}$ of system \eqref{eq_systemOpenLoop}. A parameter set of the initial condition $f_0(a),$ birth rate $k(a),$ death rate $\mu(a)$, steady-state input $D^*$ and maximum age of reproduction $A$ used throughout this work is given below
\begin{subequations}
    \label{eq_SimulationParameterAndInitialCondition}
    \begin{align}
    f_0(a) = 8 - 3a + \phi_2(a), \quad    \mu(a) = \frac{1}{20-5a}, \quad k(a) = a, \quad p(a) = 1 + \frac{a^2}{10},\label{eq_SimulationParameterSet}\\
    \Dstar \approx 0.48,\quad A=2, \quad  \phi_2(a) = \sin(\omega_1 a)\exp[\sigma_1 a]\frac{\fstar(a)}{\fstar(0)},\quad \omega_1 \approx 3.82, \quad \sigma_1 \approx 0.91. \label{eq_SimulationInitialCondition}
    \end{align}
\end{subequations}
As shown in a representative simulation example in Figure \ref{fig_ActuatorDynamics_staticController_AllFinal}, controller \eqref{eq_controllerBackstepping} achieves convergence to the desired output $\ystar = y_\mathrm{des}$. However, the controller stabilizes the equilibrium profile without restricting the dilution rate $D(t)$ to physically meaningful values.}

\section{Closed-Loop Stability for All Positive Initial Population Densities}\label{sec_BacksteppingStabilityProof}

For the stability proof, we invoke the transformation $\varPi: \mathcal{F}\to \realNumbers \times C^0([-A,0]),\ f[t] \mapsto (\eta(t),\psi_t)$ mapping the age-profile $f[t]$ to its unstable ``malthusian'' and ``asymptotic'' mode
\begin{subequations}
    \label{eq_transformation}
    \begin{align}
        \eta (t)&= \ln\Pi(f[t])\label{eta}\\
        \psi(t-a)&=\frac{f(a,t)}{\fstar(a)\Pi(f[t])}-1\label{eq_psi}
    \end{align}
\end{subequations}
where
\begin{align}
    \Pi(f) &= \frac{\int_0^A  \pi(a) f(a)\diff a}{\int_0^A ak(a)\fstar(a)\diff a},\quad \pi(a)= \int_a^A k(s)\exp[\int_s^a\left(\mu(\alpha)+\Dstar\right)\diff \alpha ]\diff s, \quad a\in[0,A].\label{eq_adjointEFpi}
\end{align} %with the notation $f[t]:[0,A]\to \realNumbers,~a\mapsto f[t](a)=f(a,t)$. 
Equation \eqref{eq_psi} has been proven to be a valid transformation \cite{KK17Chemostat}, i.e., the right hand side is indeed a function of $(t-a)$. Transforming the open-loop system \eqref{eq_systemOpenLoop} in light of \eqref{eq_transformation}, one receives
\begin{subequations}
    \label{eq_OpenLoopTransformed}
    \begin{align}
        \dot \eta (t) &= \Dstar - D(t) \label{eq_OLetaDynamics}\\ 
        \dot D (t) &= u (t) \label{eq_OLinput}\\
        \psi(t) &= \int_0^A \tilde k(a) \psi(t-a)\, \diff a ~.\label{eq_OLIDEdynamics}
    \end{align}
\end{subequations}

After transforming the {\smash{closed-loop}} system consisting of  \eqref{eq_systemOpenLoop} and \eqref{eq_controllerBackstepping}, we arrive at the dynamics
\begin{subequations}
    \label{eq_closedLoopTransformed}
    \begin{alignat}{4}
        \dot \eta (t) &= -\knom\eta (t) -\delta (t) - \knom v(\psi_t)\label{eq_etaDynamics}\\ 
        \dot \delta (t) &= - \kBackstep \delta (t)\label{eq_deltaDynamics}\\
        \psi(t) &= \int_0^A \tilde k(a) \psi(t-a)\ \diff a  \label{eq_IDEdynamics}
    \end{alignat}
\end{subequations}
where 
\begin{eqnarray}
    \tilde k(a) \coloneqq  k(a) \frac{\fstar(a)}{\fstar(0)}
    , \qquad
    g(a) \coloneqq   \frac{p(a)\fstar(a)}{\int_0^A p(a)\fstar(a)\diff a}, \qquad v(\psi_t) \coloneqq \ln \left (1+\int_0^Ag(a)\psi(t-a)\diff a \right ),
    \label{ggg}
\end{eqnarray}
 and the equilibrium $(\eta , \delta, \psi) = 0$ is to be proven asymptotically stable. This can be achieved using a member of the family of Control Lyapunov Functionals provided in Reference \citenumns{KK17Chemostat}. For this, we have to make a technical assumption:
\begin{assumption}[Technical assumption on birth kernel]
    There exists a constant $\lambda>0$ such that $\int_0^A|\tilde k(a) -r\lambda \int _a^A\tilde k(s) \diff s|\diff a<1$ where $r^{-1} = \int_0^A \tilde k (a) a \diff a$.  Let $\sigma > 0$ be a sufficiently small constant that satisfies the inequality 
    \begin{eqnarray}
    \int_0^A\left|\tilde k(a) -r\lambda \int _a^A\tilde k(s) \diff s\right| \exp[\sigma a]\diff a<1.
    \label{doy}
\end{eqnarray}
 \label{ass_technicalbirthKernel}
\end{assumption}
Assumption 1 is a mild technical assumption, since it is satisfied by arbitrary mortality rate $\mu$ for every birth kernel $k$ that has a finite number of zeros on $[0,A]$. From a modeling standpoint, this means that individuals of ages zero to $A$ are reproducing at \emph{almost every} age (see Reference \citenumns{KK17Chemostat}).

\begin{remark}
     The state $\psi$ of the internal dynamics \eqref{eq_OLIDEdynamics} is restricted (see Reference \citenumns{KK17Chemostat}) to the set
     \begin{align}
         \mathcal S = \left\{ \psi \in C^0([-A,0];(-1,\infty)): P(\psi) = 0 \land \psi(0) = \int_0^A \tilde k(a)\psi(-a)\diff a \right\}
     \end{align}
     where
     \begin{align}
         P(\psi) = \int_0^A\psi(-a)\int_a^A\tilde k(s) \diff s\diff a\left(\int_0^Aa\tilde k(a)\diff a\right)^{-1}.
     \end{align}
\end{remark}

\begin{theorem}[Lyapunov stability of backstepping controller with \emph{unrestricted dilution}] \label{thm_BacksteppingLyapunov1}
 
    % {\color{magenta}
    Let Assumption \ref{ass_technicalbirthKernel} hold. Then for every $k_{1}, k_{2} >0$, there exists a function $\alpha_1 \in \mathcal{K}_\infty$ such that for every $(f_{0},D_{0}) \in \mathcal{F} \times \realNumbers$ the unique solution $(f[t],D(t)) \in \mathcal{F} \times \realNumbers$ of the closed-loop system \eqref{eq_systemOpenLoop} with \eqref{eq_controllerBackstepping} and initial condition $(f[0],D(0))=(f_{0},D_{0})$ exists for all $t \geq 0$ and satisfies the following stability estimate for all $t \geq 0$
    \begin{align}
        R_{1}(f[t],D(t)) \leq \exp[-\frac{\sigma_1}{2} t]\;\alpha_1 \left(R_{1}(f_{0},D_{0})\right)  \label{eq_expBoundOriginalStateDtilde}
    \end{align}
    where $\sigma_1 = \min\left(\frac{\KarafControlGain{1}}{2},\KarafControlGain{2},\sigma\right) > 0$ and
     \begin{align}
     R_{1}(f,D):= \max_{a\in [0,A]}\left| \ln \frac{f(a)}{\fstar(a)} \right| + \left| D-\Dstar\right|
    \label{eq_measureR1(f,F)}
    \end{align}
    for all $(f,D) \in \mathcal{F} \times \realNumbers$.
    % }
\end{theorem}

\begin{proof}[Proof of Theorem \ref{thm_BacksteppingLyapunov1}]
    Lemma 4.1 of Reference \citenumns{KK17Chemostat} provides existence and uniqueness of the solution $\psi_t \in \mathcal S$ and also
    \begin{equation}
        \inf_{t\geq -A} \psi(t) \geq \min_{t\in [-A,0]} \psi(t)> -1, \quad \forall t\geq 0. \label{eq_IDEstatebound}
    \end{equation}
    Notice that the second inequality is given by $\psi_0 \in \mathcal S$ and is guaranteed for all physical initial conditions $f_0 \in\mathcal{F}$. Indeed, the initial condition of \eqref{eq_IDEdynamics} is given by
    \begin{equation}
        \psi_0(a) = \psi(-a) = \frac{f_0(a)}{\fstar (a) \Pi (f_0)} -1 > -1, \quad \forall  a \in [0,A] \label{eq_functionalPsiV}
    \end{equation}
    where $f_0(a)$, $\fstar (a)$ and $\Pi(f_0)$ are positive on the domain $[0,A]$. Since $g(a) \geq 0$, $\int_0^A g(a)\diff a = 1$ according to \eqref{ggg} and inequality \eqref{eq_IDEstatebound} holds, the map 
    \begin{equation}\label{eq-vdef}
        %t \mapsto 
        v(\psi_t) \coloneqq \ln \left (1+\int_0^Ag(a)\psi(t-a)\diff a \right ), \quad \textrm{for} \quad t \geq 0
    \end{equation}
    is well-defined and  continuous. Thus, the ODE subsystem of the closed-loop system \eqref{eq_closedLoopTransformed}, namely, equations \eqref{eq_etaDynamics}, \eqref{eq_deltaDynamics} locally admits a unique solution.
%We next show \eqref{eq_DiniBoundLyap}. 
    Consider first the Lyapunov function
    \begin{equation}
        U_1(\eta,\delta) = \frac{1}{2}\left(\eta^2+\LyapCoeff{1}\delta^2\right), \quad \forall(\eta,\delta) \in \realNumbers^2
    \end{equation}
    where $\LyapCoeff{1}>0$ is a constant to be chosen. The time derivative of $U_1(\eta(t),\delta(t))$ along the solutions of \eqref{eq_etaDynamics}, \eqref{eq_deltaDynamics} can be upper bounded for all times $t \geq 0$ for which the solution $(\eta(t),\delta(t))$ exists as follows:
    \begin{align}
        \frac{\diff}{\diff t}U_1(\eta(t),\delta(t)) &\leq -\frac{\knom}{4} \eta(t)^2 - \left(\LyapCoeff{1}\kBackstep-\frac{1}{\knom}\right) \delta(t)^2 + \frac{\knom}{2} v(\psi_t)^2. \label{eq_BoundODELyap}
    \end{align}
    Inequality \eqref{eq_BoundODELyap} was derived using the inequalities $-\eta(t) \delta(t) \leq \frac{\knom}{4} \eta(t)^2 + \frac{1}{\knom} \delta(t)^2$ and $- \eta(t) v(\psi_t) \leq \frac{1}{2} v(\psi_t)^2 + \frac{1}{2} \eta(t)^2$.
Next, as in the proof of Reference \citenumns[Lemma 5.1]{KK17Chemostat}, we use Corollary 4.6 of the latter with $C(\psi_t) = (1+\min (0,\min_{a\in [0,A] }\psi_t(a)))^{-2}$ and $b(s) = \frac{1}{2} s^2$ under Assumption \ref{ass_technicalbirthKernel} with sufficiently small parameter $\sigma >0$ to obtain that
        the functional
    \begin{align}\label{eq-G}
        G(\psi_t) \coloneqq \frac{\max_{a\in [0,A]} |\psi(t-a)|\exp[-\sigma a]}{1+\min(0,\min_{a\in [0,A]}\psi(t-a))}
    \end{align}
    satisfies
    \begin{equation}
        D^+ G(\psi_t)^2 \leq -2 \sigma  G(\psi_t)^2, \quad \forall t \geq 0\label{eq_BoundPsiLyap}
    \end{equation}
    along solutions $\psi_t$ of the IDE \eqref{eq_IDEdynamics}.
    % where $D^+$ denotes the upper right Dini-derivative
    % \begin{equation}
    %     D^+ f(t) = \limsup_{h\to 0^+}\left\{h^{-1}(f(t+h)-f(t))\right\}
    % \end{equation}
    % with $f\in C^0(\realNumbers_{\geq 0}; \realNumbers).$
    Now, Reference \citenumns[(A.43)]{KK17Chemostat} also provides that 
    \begin{equation}
            |v(\psi_t)|\leq G(\psi_t)\exp[\sigma A], \quad \forall t \geq 0 .\label{eq-vleqG}
    \end{equation}
 Given this fact, the bounds \eqref{eq_BoundODELyap} and \eqref{eq_BoundPsiLyap} and defining the Lyapunov functional
    \begin{eqnarray}
        V_1(\eta,\delta,\psi) &=& U_1(\eta,\delta) + \frac{1}{2}\LyapCoeff{2} G(\psi)^2, 
    % \end{equation}
    %     \begin{align}
%        & V_1(\eta,\delta,\psi) 
        \nonumber\\
        &=& \frac{1}{2}\left(\eta^2 + \LyapCoeff{1}\delta^2 + \LyapCoeff{2} G(\psi)^2\right), \quad \forall (\eta,\delta,\psi)\in \realNumbers \times \realNumbers \times \mathcal{S}\label{eq_LyapV1Definition}
    \end{eqnarray}
    with 
        $\LyapCoeff{1} = \frac{2}{\kBackstep\knom}$, $\LyapCoeff{2} = \frac{\knom}{\sigma}\exp[2\sigma A]$, 
        we derive the differential inequality
    \begin{equation}
        D^+ V_1(\eta(t),\delta(t),\psi_t) \leq -\sigma_1 V_1(\eta(t),\delta(t),\psi_t),\quad \forall t\geq 0 \label{eq_DiniBoundLyap}
    \end{equation}
    %     \begin{equation}
    %     D^+ V_1(\eta(t),\delta(t),\psi_t) \leq -\sigma_1 V_1(\eta(t),\delta(t),\psi_t),~\text{ for all }t\geq 0 \label{eq_DiniBoundLyap}
    % \end{equation}
    % where $\sigma_1 = \min\left(\frac{1}{2}\knom,\kBackstep,\sigma\right) > 0$.
     % estimate \eqref{eq_DiniBoundLyap}, where one chooses 
    % $\LyapCoeff{1} = \frac{2}{\kBackstep\knom}$, $\LyapCoeff{2} = \frac{\knom}{\sigma}\exp[2\sigma A]$ sufficiently large and $\gamma_1 = 1$, $\gamma_2 = \frac{1}{2} \knom$ sufficiently small 
    % and receives $\sigma_1 = \min\left(\frac{1}{2}\knom,\kBackstep,\sigma\right)$. 
        where $\sigma_1 = \min\left(\frac{\KarafControlGain{1}}{2},\KarafControlGain{2},\sigma\right) > 0$. Since estimate \eqref{eq_DiniBoundLyap} implies that all solutions stay bounded for all times for which they exist, we can conclude existence for all times. The differential inequality \eqref{eq_DiniBoundLyap} in conjunction with Lemma 2.12 on pages 77-78 in Reference \citenumns{karafyllisJiang2011stability} implies the following estimate for all $t \geq 0$:
    \begin{align}
         V_1(\eta(t),\delta(t),\psi_t) \leq \exp[-\sigma_1 t] V_1(\eta(0),\delta(0),\psi_{0}).
        \label{doula}
    \end{align}

{
% SHOWING BOUND IN $R(f,\tilde D)$: 
We next show that the estimate \eqref{eq_expBoundOriginalStateDtilde} holds. Recall from \eqref{eq_defintionDelta} that 
\begin{align}
    \delta(t) &= \tilde D(t) - \KarafControlGain{1} \left( \eta(t) + v(\psi_t) \right), 
\quad  t \geq 0\label{eq_deltaInverseTransform}
\end{align}
where
\begin{align}
    \tilde D(t) = D(t) - \Dstar, \quad t \geq 0.
\end{align}
Define the positive definite quadratic function $U_2: \realNumbers^3 \to \realNumbers$
\begin{align}
    U_2(\eta,\tilde D, v) = \eta^2 + \frac{\LyapCoeff{2}}{2} \exp[-2\sigma A] v^2 + \left( \tilde D - \KarafControlGain{1}\eta - \KarafControlGain{1}v \right)^2
\end{align}
where $v$ is given by \eqref{eq-vdef}. 
Clearly, there exist constants $\underline{c},\overline{c} >0$ such that for all $ (\eta,\tilde D, v)\in \realNumbers^3$ the following inequalities hold:
\begin{align}
    \underline{c}\left(\eta^2 + \tilde D^2 + v^2 \right) \leq U_2(\eta,\tilde D, v) \leq \overline{c}\left(\eta^2 + \tilde D^2 + v^2 \right). \label{eq_upperAndLowerBoundU2}
\end{align}
 Consequently, \eqref{eq_upperAndLowerBoundU2} in conjunction with \eqref{eq_LyapV1Definition}, \eqref{eq_deltaInverseTransform}, \eqref{eq_upperAndLowerBoundU2}, \eqref{eq-vleqG} implies the following estimate:
\begin{align}
    V_1(\eta(t),\delta(t),\psi_t) &= \frac{1}{2}\left(U_2(\eta(t),\tilde D(t), v(\psi_t)) + \LyapCoeff{2} G(\psi_t)^2 - \frac{\LyapCoeff{2}}{2} \exp[-2\sigma A] v(\psi_t)^2\right)\\
    &\geq \frac{1}{2}\left(\underline{c}\left(\eta(t)^2 + \tilde D(t)^2 + v(\psi_t)^2 \right) + \frac{\LyapCoeff{2}}{2} G(\psi_t)^2 \right), \quad t \geq 0.\label{eq_lyapDtildebound}
\end{align}
Using \eqref{eq_lyapDtildebound}, we conclude that there exists an appropriate constant $\tilde c>0$ (independent of the solution), for which the following estimate holds:
\begin{align}
    |\eta(t)| + \exp[\sigma A] G(\psi_t) + \left|\tilde D(t)\right| &\leq \tilde c\sqrt{ V_1(\eta(t),\delta(t),\psi_t)}, \quad t \geq 0 .\label{eq_lowerBoundV_1}
\end{align}
From Reference \citenumns[(5.24)]{KK17Chemostat},
    \begin{align}
        \max_{a\in [0,A]}\left| \ln \frac{f[t](a)}{\fstar(a)} \right| \leq |\eta(t)| + \exp[\sigma A] G(\psi_t), \quad t \geq 0,\label{eq_KK17_EtaPsiBoundToF_upper}
    \end{align}
    i.e., we arrive at the left-hand side of \eqref{eq_expBoundOriginalStateDtilde} after combining \eqref{eq_lowerBoundV_1} and \eqref{eq_KK17_EtaPsiBoundToF_upper}. To complete the proof of \eqref{eq_expBoundOriginalStateDtilde}, we need to show that there exists $ \tilde\alpha \in \mathcal{K}_\infty$ such that 
    \begin{align}
        V_1(\eta(0),\delta(0),\psi_0) \leq \tilde\alpha\left(\max_{a\in [0,A]}\left| \ln \frac{f_0(a)}{\fstar(a)} \right| + \left|\tilde D(0)\right|\right). \label{eq_initialLyapboundClassK}
    \end{align}
  We first note that Reference \citenumns[(5.27)]{KK17Chemostat} provides the existence of $\tilde\alpha_1,~\tilde\alpha_2 \in \mathcal{K}_\infty$
    \begin{align}
        |\eta_0| \leq \tilde\alpha_1\left(\max_{a\in [0,A]}\left| \ln \frac{f_0(a)}{\fstar(a)} \right| \right),\quad G(\psi_0) \leq \tilde\alpha_2 \left(\max_{a\in [0,A]}\left| \ln \frac{f_0(a)}{\fstar(a)} \right|\right), \label{eq_KK17technicalbounds}
    \end{align}
and, using \eqref{eq_upperAndLowerBoundU2}, we upper bound
\begin{align}
    V_1(\eta,\delta,\psi) \leq \frac{1}{2}\left(\overline{c}\left(\eta^2  + \tilde D(\eta,\delta,\psi)^2 \right) + \tilde b G(\psi)^2 \right),\quad (\eta,\delta,\psi) \in \realNumbers\times \realNumbers \times \mathcal{S}\label{eq_boundingforalphatilde}
\end{align}
for an appropriate constant $\tilde b > 0$. The inequalities \eqref{eq_KK17technicalbounds} and \eqref{eq_boundingforalphatilde} allow us to conclude that there exists appropriate $\tilde \alpha$ (independent of the solution) for which \eqref{eq_initialLyapboundClassK} holds. Combining \eqref{doula}, \eqref{eq_lowerBoundV_1}, \eqref{eq_KK17_EtaPsiBoundToF_upper} and \eqref{eq_initialLyapboundClassK}, we obtain estimate \eqref{eq_expBoundOriginalStateDtilde}.
}
\end{proof}
\begin{remark}
    Notice that the exponential convergence rate $\sigma_1$ of \eqref{eq_DiniBoundLyap} is the minimum of three constants: the nominal control gain $\knom$, the backstepping convergence rate $\kBackstep$ and the convergence rate of the internal dynamics $\sigma$.
\end{remark}

\section{Easier-to-Implement Backstepping Controller: Output Feedback without Sensor Model} \label{sec_RelaxedBackstepping}

\begin{figure}[t]
    \centering
    \includegraphics[width=\textwidth]{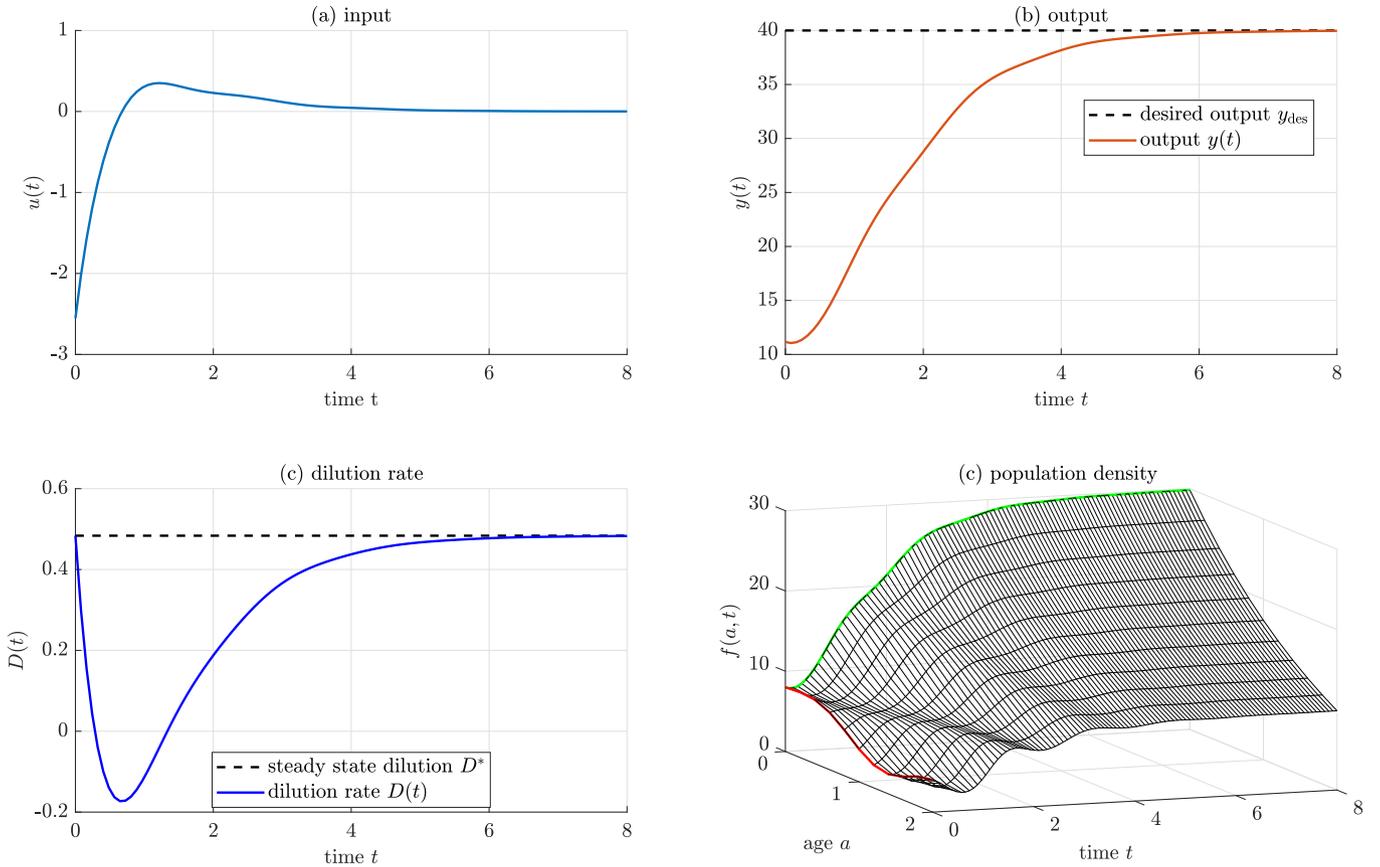}
    \caption{Simulation results of system \eqref{eq_systemOpenLoop} with parameters and initial condition \eqref{eq_SimulationParameterAndInitialCondition} and controller \eqref{eq_controllerBackstepping} but relaxed canceling terms \eqref{eq_canceling_Relaxed} with gains $\KarafControlGain{1} = 1$, $\KarafControlGain{2} = 2$. Convergence to the desired output $\ystar = y_\mathrm{des}$ is achieved. The controller employs a non-physical negative dilution $D(t)$ in the transient to achieve stabilization globally (for all positive population densities). A more complicated backstepping re-design that employs global stabilization using only positive dilution is presented in Section \ref{sec_Safety}.}
    \label{fig_ActuatorDynamics_staticControllerRelaxed_AllFinal}
\end{figure}

Recall the canceling controller \eqref{eq_cancelinggeneral}
\begin{equation}
     \uCancel(t) = \knom\left( -D(t)+\Dstar + v_1(\psi_t)\right)\label{eq_cancelinggeneral_recall}\\
\end{equation}
where the tedious term 
\begin{multline}
    v_1(\psi_t) = -\Dstar + \frac{1}{\int_0^A p(a) \fstar (a) (1+\psi(t-a))\diff a} \\ 
    \left [ p(A) \fstar (A) (1+\psi(t-A)) - p(0) \fstar (0) (1+\psi(t)) - \int_0^A \tilde p(a) \fstar (a) (1+\psi(t-a))\diff a \right] \label{eq_functionalPsi}
\end{multline}
is a functional of the internal state variable $\psi$ and thus decays exponentially to $v_1(0) = 0$. We relax the canceling control from exact cancellation to now only canceling the steady-state value 
\begin{equation}
    \uCancel(t) = \knom\left(-D(t) + \Dstar\right) \label{eq_canceling_Relaxed}
\end{equation}
drastically reducing the measurement requirements. %Moreover, this simplification will carry over to the safety filter as well. 
This approach achieves convergence $y(t)\to \ystar$ to the desired output, which is illustrated by means of a representative simulation example in Figure \ref{fig_ActuatorDynamics_staticControllerRelaxed_AllFinal}. 

Finding the closed-loop dynamics of \eqref{eq_systemOpenLoop}, \eqref{eq_controllerBackstepping} and \eqref{eq_canceling_Relaxed}
\begin{subequations}
    \label{eq_closedLoopSuperlateBotchingTransformed}
    \begin{alignat}{4}
        \dot \eta (t) &= -\knom\eta (t)  -\delta (t) - \knom v(\psi_t)\label{eq_etaDynamicsSuperlateBotching}\\ 
        \dot \delta (t) &= -\kBackstep\delta(t) -\knom v_1(\psi_t)\label{eq_deltaDynamicsSuperlateBotching}\\
        \psi(t) &= \int_0^A \tilde k(a) \psi(t-a)\, \diff a,  \label{eq_IDEdynamicsSuperlateBotching}
    \end{alignat}
\end{subequations}
we notice that, as a result of simplifying the controller, the dynamics of $\delta$ have become more complex, and as a consequence the Lyapunov-like stability proof becomes harder to establish. 

In our proof of stability, we make use of the following result, whose proof we state in the Appendix.
\begin{lemma}\label{lemma_T_1_boundMain}
    There exists a constant $c_1>0$ such that every solution $\psi_t$ of the IDE \eqref{eq_IDEdynamicsSuperlateBotching} satisfies the following estimate for all $t\geq 0$: 
    \begin{align} \label{eq_T_1_boundMain}
        |v_1(\psi_t)| \leq \sqrt{c_1}\frac{ \max_{a\in[0,A]}|\psi(t-a)|}{1+\min_{a\in [0,A]}\psi(t-a)}.%\leq \frac{c_1}{\ystar}\frac{||\psi_t||_\infty}{1+\min_{a\in [0,A]}\psi_0(a)}
    \end{align}
\end{lemma}
\begin{theorem} [Lyapunov stability of \emph{static output feedback} backstepping controller with unrestricted dilution]\label{thm_BacksteppingRelaxed2}
    %
    % {\color{magenta}
    Let Assumption \ref{ass_technicalbirthKernel} hold. Then for every $k_{1}, k_{2} >0$, there exists a function $\alpha_2 \in \mathcal{K}_\infty$ such that for every $(f_{0},D_{0}) \in \mathcal{F} \times \realNumbers$ the unique solution $(f[t],D(t)) \in \mathcal{F} \times \realNumbers$ of the closed-loop system \eqref{eq_systemOpenLoop} %, \eqref{eq_controllerBackstepping+} 
    %with \eqref{eq_canceling_Relaxed} 
    with 
\begin{equation}    \label{eq_controllerBackstepping+}
u(t) = - (\knom+\kBackstep)\left(D(t) - \Dstar\right)
+ \knom \kBackstep \ln \frac{y(t)}{\ystar} 
\end{equation}
    and initial condition $(f[0],D(0))=(f_{0},D_{0})$ exists for all $t \geq 0$ and satisfies the following stability estimate for all $t \geq 0$
    \begin{align}
        R_{1}(f[t],D(t)) \leq \exp[-\frac{\sigma_1}{2} t]\;\alpha_2 \left(R_{1}(f_{0},D_{0})\right)  \label{eq_expBoundOriginalStateDtildeRelaxedController}
    \end{align}
    where $\sigma_1 = \min \left(\frac{\KarafControlGain{1}}{2},\frac{\KarafControlGain{2}}{2},\sigma\right) > 0$ and $R_1$ is defined in \eqref{eq_measureR1(f,F)}. 
\end{theorem}

\begin{proof}[Proof of Theorem \ref{thm_BacksteppingRelaxed2}]
    Again, like in \eqref{eq_functionalPsiV}, we know that the map $t \mapsto v(\psi_t)$ is well-defined for all $t \geq 0$ {and is continuous}. Thus, the closed-loop system \eqref{eq_closedLoopSuperlateBotchingTransformed} locally admits a solution. We first show \eqref{eq_DiniBoundLyapRelaxed} by considering the function
    \begin{equation}
        U_1(\eta,\delta) = \frac{1}{2}(\eta^2+\LyapCoeff{1}\delta^2), \quad (\eta, \delta) \in \realNumbers^2,
    \end{equation}
    {where $\LyapCoeff{1}>0$ is a constant to be chosen}. The time derivative of {$U_1(\eta(t),\delta(t))$} along solutions of \eqref{eq_etaDynamicsSuperlateBotching}, \eqref{eq_deltaDynamicsSuperlateBotching} can be upper bounded for all times $t \geq 0$ for which the solution $(\eta(t),\delta(t))$ exists
    \begin{align}
        \frac{\diff}{\diff t}U_1(\eta(t),\delta(t)) &\leq -\frac{\knom}{4}\eta(t)^2 - \left(\frac{1}{2}\LyapCoeff{1}\kBackstep-\frac{1}{\knom}\right)\delta(t)^2 \notag\\
        &+ \frac{\knom}{2}v(\psi_t)^2 + \LyapCoeff{1}\frac{\knom^2}{2\kBackstep}v_1(\psi_t)^2 \label{eq_BoundODELyapRelaxedFinal}.
    \end{align}
    % \begin{align}
    %     \frac{\diff}{\diff t}U_1(\eta(t),\delta(t)) &\leq -\left(\knom\left(1-\frac{1}{2}\gamma_1\right)-\frac{1}{2}\gamma_2\right)\eta(t)^2 - \left(\LyapCoeff{1}\left(\kBackstep-\frac{1}{2}\gamma_3\right)-\frac{1}{2\gamma_2}\right)\delta(t)^2 \notag\\
    %     &+ \frac{\knom}{2\gamma_1}\exp[2\sigma A]G(\psi_t)^2 + \LyapCoeff{1}\frac{c_1^2\knom^2}{2\gamma_3 {\ystar}^2}\exp[2\sigma A]G_2(\psi_t)^2 \label{eq_BoundODELyapRelaxedFinal}.
    % \end{align}
    {The bound \eqref{eq_BoundODELyapRelaxedFinal} was achieved using the inequalities $-\eta(t)\delta(t)\leq \frac{\KarafControlGain{1}}{4}\eta(t)^2 + \frac{1}{\KarafControlGain{1}}\delta(t)^2$, $-\eta(t)v(\psi_t) \leq \frac{1}{2}\eta(t)^2 + \frac{1}{2}v(\psi_t)^2$ and $-\delta(t)v_1(\psi_t) \leq \frac{\KarafControlGain{2}}{2}\delta(t)^2 + \frac{1}{2\KarafControlGain{2}}v_1(\psi_t)^2$.}
    % using Young's inequality for the indefinite terms with the to-be-chosen constants $\gamma_i >0$, 

    Next, consider the Lyapunov functional
    \begin{subequations}
        \begin{align}
            & V_2(\eta(t),\delta(t),\psi_t) = U_1(\eta(t),\delta(t)) +  \frac{\LyapCoeff{2}}{2} G(\psi_t)^2  = \frac{1}{2} \left( \eta(t)^2 + \LyapCoeff{1}\delta(t)^2 + \LyapCoeff{2} G(\psi_t)^2 \right) \label{eq_lyapunovFunctionalRelaxedControl}
        \end{align}
        where $\LyapCoeff{2}>0$ is a constant to be chosen and $G$ is defined in \eqref{eq-G}.
        % \begin{align}
        %     &G(\psi_t) \coloneqq \frac{G_0(\psi_t)}{1+\min(0,\min_{a\in [0,A]}\psi(t-a))},\quad G_1(\psi_t) \coloneqq \frac{\max_{a\in [0,A]} |\psi(t-a)|\exp[-\sigma a]}{1+\min_{a\in [0,A]}\psi(t-a)}, \\
        %     & G_0(\psi_t) = \max_{a\in [0,A]} |\psi(t-a)|\exp[-\sigma a].
        % \end{align}
        % \begin{align}
        %    G_1(\psi_t) \coloneqq \frac{\max_{a\in [0,A]} |\psi(t-a)|\exp[- \sigma a]}{1+\min_{a\in [0,A]}\psi(t-a)}.
        % \end{align}
    \end{subequations}

    Using Lemma \ref{lemma_T_1_boundMain} and definition \eqref{eq-G}, it holds that
    \begin{align}\label{eq_v2ConnectionToG_2}
        |v_1(\psi_t)|\leq \exp[\sigma A]G(\psi_t), \quad \forall t\geq 0.
    \end{align}
    
    Using the above estimate, \eqref{eq_BoundODELyapRelaxedFinal}, \eqref{eq_BoundPsiLyap} and defining $\LyapCoeff{1} = \frac{4}{\kBackstep\knom}$, $\LyapCoeff{2} = \left( \frac{\knom}{\sigma} + \LyapCoeff{1}\frac{c_1\knom^2}{\sigma \kBackstep } \right) \exp[2\sigma A]$, we obtain the differential inequality
    \begin{equation}
        D^+ V_2(\eta(t),\delta(t),\psi_t) \leq -\sigma_1 V_2(\eta(t),\delta(t),\psi_t),\quad \forall t\geq 0 \label{eq_DiniBoundLyapRelaxed}
    \end{equation}
    where $\sigma_1 = \min\left (\frac{1}{2}\knom,\frac{1}{2}\kBackstep,\sigma\right ) > 0$. Since estimate \eqref{eq_DiniBoundLyapRelaxed} implies that all solutions stay bounded for all times for which they exist, we can conclude existence for all times.
    
    The rest of the proof is exactly the same with the proof of Theorem \ref{thm_BacksteppingLyapunov1}.
\end{proof}

\section{Global Stabilization with State Constraints}

\label{sec_Safety}

The proposed controller \eqref{eq_controllerBackstepping} unfortunately does not ensure that the dilution rate $D(t)$ remains positive. 
%, in that there are sets of initial conditions and desired outputs such that we receive negative Dilution rates, which are non-physical, as can be seen 
This can be observed in an example simulation in Figure \ref{fig_ActuatorDynamics_staticController_AllFinal}. {Whenever the signed control error $y(t)-\ystar$ takes large negative values, the controller tries to ``add population.''}

A classical approach to guarantee positivity of the dilution rate is to introduce a control barrier function (CBF) 
\begin{equation}
    h(\eta, D, \psi) = D ,
    %\overset{!}{>} 0~.
\end{equation} 
which imposes a positivity constraint on the state $D$, and a safety filter of the nominal control law $u_0= \uCancel + \uStabilize$, so that the overall control is given by
\begin{equation}
    u = u_0 + u_\mathrm p = \uCancel + \uStabilize + u_\mathrm p 
    \label{eq_SafetyAugmentedControl}
\end{equation}
where $u_\mathrm p$ is the safety filter/override component of the controller that ``penalizes'' the negativity of the dilution and overrides the nominal feedback $u_0$ so that the dilution is kept positive. 

The safety filter for a system that has actuator dynamics $\dot D =u$ and CBF  $h = D$ is particularly simple and given in Reference \citenumns{10130642}
\begin{equation}
    u_\mathrm p = \max \{0,-u_0 - \kSafe D\}, \quad \kSafe>0, 
    \label{eq_safetyFilter}
\end{equation}
which ensures that 
% Applying the controller \eqref{eq_SafetyAugmentedControl} and \eqref{eq_safetyFilter} to the system \eqref{eq_systemOpenLoop} ensures
\begin{align}
    \dot h(t) &\geq - \kSafe h(t),
    %,\quad\forall  x \in \realNumbers\times\realNumbers\times\mathcal{S}
\end{align}
namely, that $\dot D(t) \geq - \kSafe D(t)$, i.e., $D(t) \geq D(0) \exp[-\kSafe t] >0$ for $D(0)>0$. 

\begin{figure}[t]
    \centering
    \includegraphics[width=\textwidth]{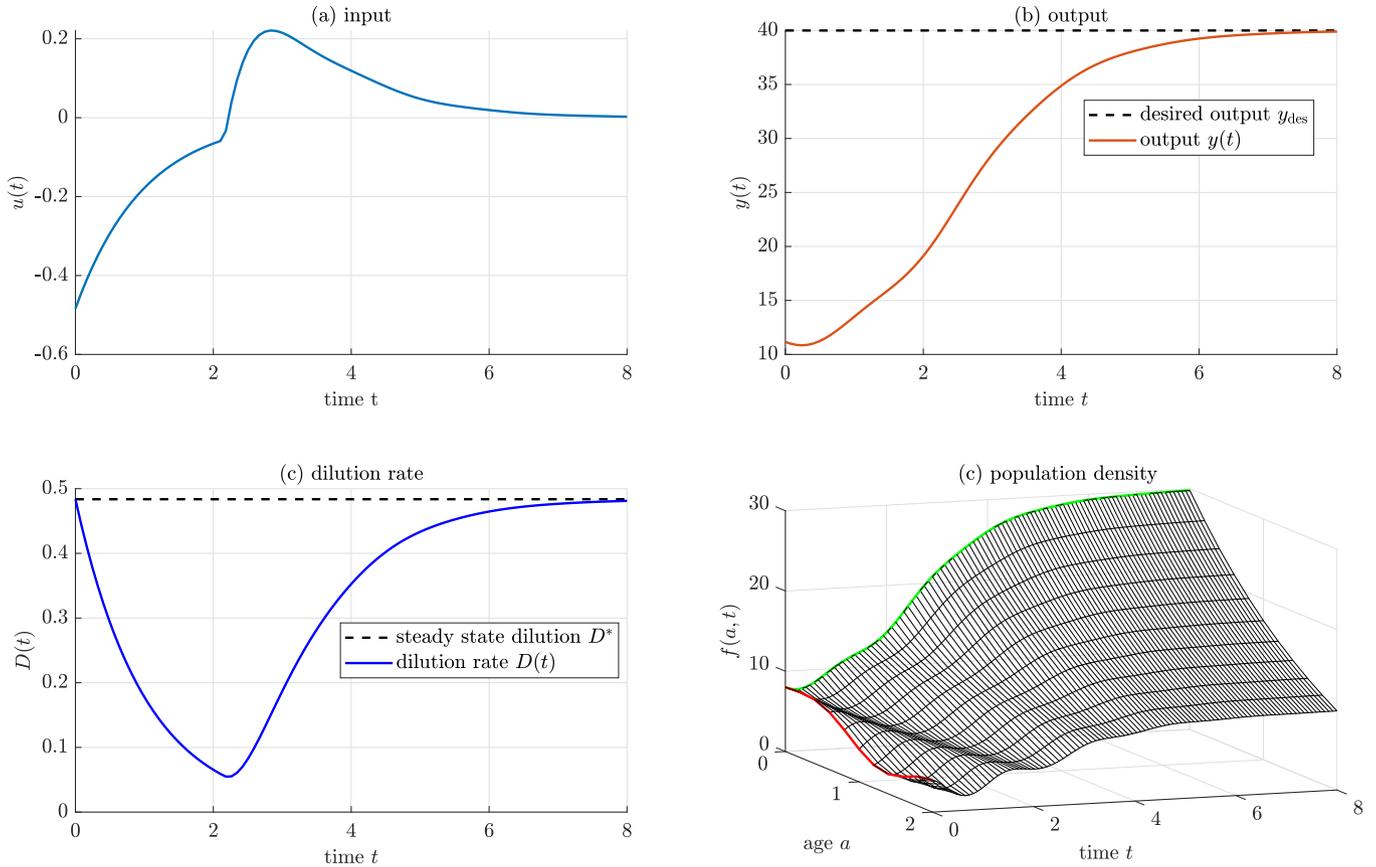}
    \caption{Simulation results of system \eqref{eq_systemOpenLoop} with parameters and initial condition \eqref{eq_SimulationParameterAndInitialCondition} and controller \eqref{eq_SafetyAugmentedControl} and \eqref{eq_safetyFilter} with gains $\KarafControlGain{1} = 1$, $\KarafControlGain{2} = 2$, $\KarafControlGain{3} = 1$. Negative dilution rates are prevented and convergence to the desired output $\ystar = y_\mathrm{des}$ is achieved.}
    % Simulation results of system \eqref{eq_systemOpenLoop} with controller \eqref{eq_SafetyAugmentedControl} and \eqref{eq_safetyFilter} and parameters $\mu = 0.2$, $p = 1$, $k(a) = 2a(A-a)$ and $A = 2$. Negative dilution rates are prevented and stabilization is achieved.
    \label{fig_ActuatorDynamics_staticControllerSafety_AllFinal}
\end{figure}

The representative example simulation in Figure \ref{fig_ActuatorDynamics_staticControllerSafety_AllFinal} shows that the safety filter ensures the positivity of $D(t)$ without impacting the convergence to the desired equilibrium. However, proving that the region of attraction of an equilibrium contains the entire safe set in the presence of a safety filter is impossible in general and has eluded us with this particular system as well. 

For this reason, we investigate the possibility of ensuring the positivity of $D(t)$ with a stabilizing controller different than the nominal control \eqref{eq_controllerBackstepping}. We next study system \eqref{eq_systemOpenLoop} with state space $\mathcal{F} \times (\Dmin,\Dmax)$, where $0 \leq \Dmin<\Dstar<\Dmax$ are positive constants. The values of the constants $\Dmin, \Dmax$ are determined by the technical characteristics of the bioreactor but from a mathematical point of view the values of the constants $\Dmin, \Dmax$ are considered to be arbitrary given constants.  

Consider the diffeomorphism $\varPhi:\realNumbers \to (\Dmin,\Dmax)$ defined by
\begin{align}
    D = \varPhi(\zeta) = \Dmin + \frac{\KarafParA \exp[\zeta]}{\KarafParB+\exp[\zeta]}
    \label{diffeo1}
\end{align}
with inverse $\varPhi^{-1}:(\Dmin,\Dmax) \to \realNumbers$ 
\begin{align}
    \zeta = \varPhi^{-1}(D) = \ln \left (\frac{\KarafParB(D-\Dmin)}{\Dmax-D} \right )
    \label{diffeo2}
\end{align}
where
\begin{align}
    \KarafParA &= \Dmax - \Dmin > 0, ~~\KarafParB = \frac{\Dmax-\Dstar}{\Dstar-\Dmin} > 0.
    \label{diffeo3}
\end{align}
Notice that by virtue of \eqref{diffeo1}, \eqref{diffeo2} and \eqref{diffeo3} it holds that $\varPhi^{-1}(\Dstar)=0$ and $\varPhi(0)=\Dstar$.

Let $k_{1}, k_{2}, k_{3}>0$ be given arbitrary positive constants (the controller gains) and consider the static output feedback law
\begin{align}
    u(t) = \frac{(D(t)-\Dmin)(\Dmax-D(t))}{\Dmax - \Dmin}\left [(\KarafControlGain{1} + \KarafControlGain{2}) (\Dstar-D(t)) - \KarafControlGain{3}\left ( \varPhi^{-1}(D(t))- \KarafControlGain{2} \ln (y(t)/\ystar) \right) \right]. \label{eq_IAS1_controller}
\end{align}
Our main result in this section is stated next.

\begin{figure}[t]
    \centering
    \includegraphics[width=\textwidth]{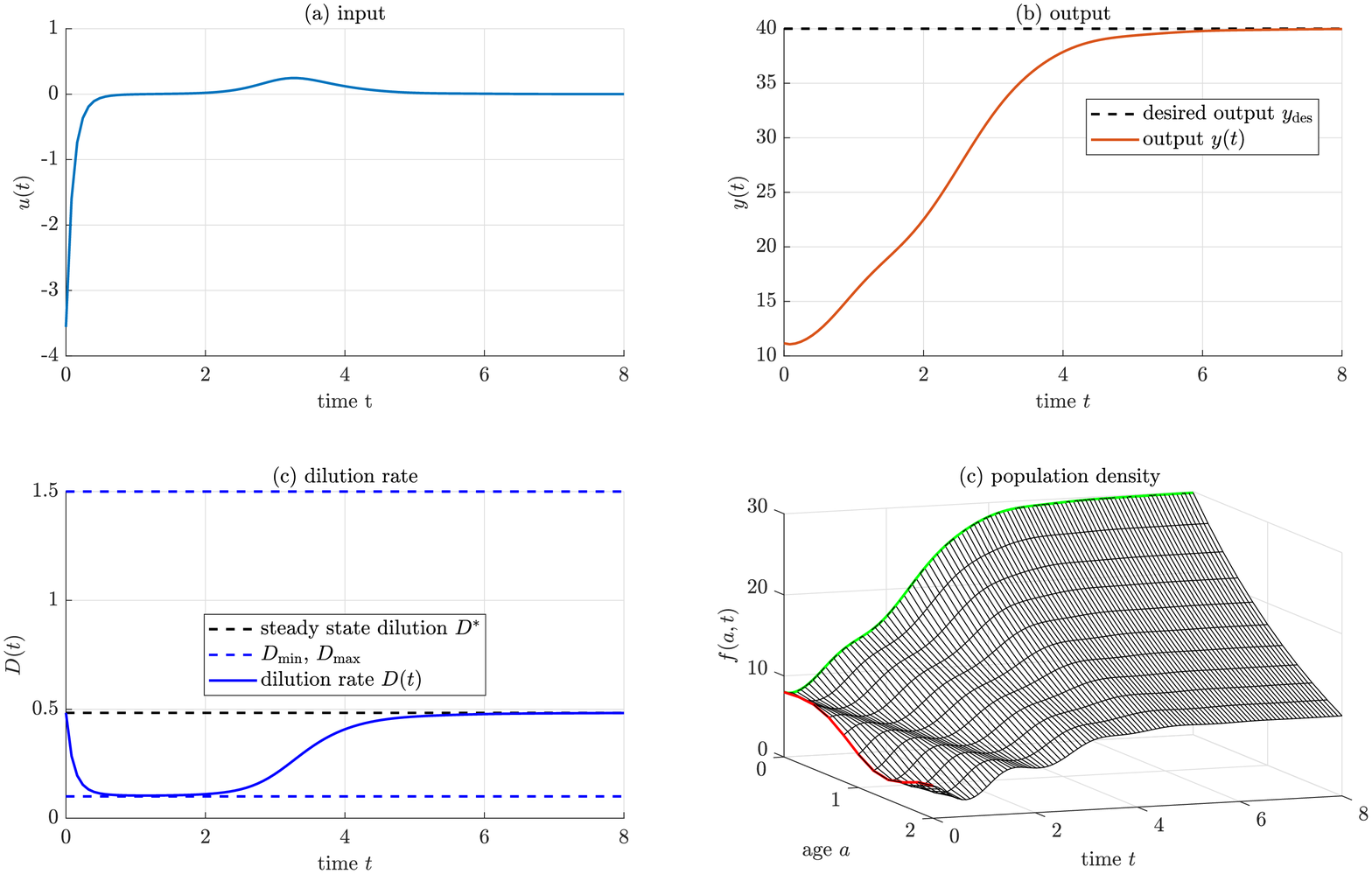}
    \caption{Simulation results of system \eqref{eq_systemOpenLoop} with parameters and initial condition \eqref{eq_SimulationParameterAndInitialCondition} and controller \eqref{eq_IAS1_controller} with gains $\KarafControlGain{1} = 1$, $\KarafControlGain{2} = 10$, $\KarafControlGain{3} = 1$ and dilution interval $(\Dmin,\Dmax) = (D_\mathrm{min},D_\mathrm{max}) = (0.1,1.5)$. The controller constrains the dilution rate $D(t) \in (\Dmin,\Dmax)$ and convergence to the desired output $\ystar = y_\mathrm{des}$ is achieved.}
    \label{fig_simulationActDynKarafyllis}
\end{figure}

\begin{theorem} [Global asymptotic stabilization with dilution constrained to a finite positive interval] \label{thm_ConstrainedDilutionGlobalAsymStability}
    Let Assumption \ref{ass_technicalbirthKernel} hold. Then for every $k_{1}, k_{2}, k_{3}>0$, there exists a non-increasing function $\varphi : \realNumbers_+ \to (0,+\infty) $ and a function $\alpha \in \mathcal{K}_\infty$ such that for every $(f_{0},D_{0}) \in \mathcal{F} \times (\Dmin,\Dmax)$ the unique solution $(f[t],D(t)) \in \mathcal{F} \times (\Dmin,\Dmax)$ of the closed-loop system \eqref{eq_systemOpenLoop} with \eqref{eq_IAS1_controller} and initial condition $(f[0],D(0))=(f_{0},D_{0})$ exists for all $t \geq 0$ and satisfies the following stability estimate for all $t \geq 0$
    \begin{align}
        R_{2}(f[t],D(t)) \leq \exp[- \varphi \left( R_{2}(f_{0},D_{0}) \right) t] \alpha \left( R_{2}(f_{0},D_{0}) \right) 
        \label{eq_IAS2_DecayBound}
    \end{align}
    where
     \begin{align}
     R_{2}(f,D):=\max_{a\in [0,A]}\left| \ln \frac{f(a)}{\fstar(a)} \right| + |\varPhi^{-1}(D)|
    \label{eq_IAS2_DecayBound22}
    \end{align}
    for all $(f,D) \in \mathcal{F} \times (\Dmin,\Dmax)$.
\end{theorem}

\emph{{Discussion of Theorem \ref{thm_ConstrainedDilutionGlobalAsymStability}:}} The stability estimate \eqref{eq_IAS2_DecayBound} shows Global Asymptotic Stability of the equilibrium point $(\fstar,\Dstar)$ for the closed-loop system \eqref{eq_systemOpenLoop} with \eqref{eq_IAS1_controller}. However, there are important differences with the corresponding stability estimate of Theorem \ref{thm_BacksteppingRelaxed2}. 

\begin{itemize}
  \item The equilibrium point $(\fstar,\Dstar)$ for system \eqref{eq_systemOpenLoop} with \eqref{eq_IAS1_controller} is stabilized in a special measure: the measure $R_{2}(f,D)=\max_{a\in [0,A]}\left| \ln (f(a)/\fstar(a)) \right| + |\varPhi^{-1}(D)|$ and not the measure $R_{1}(f,D)=\max_{a\in [0,A]}\left| \ln (f(a)/\fstar(a)) \right| + |D-\Dstar|$ that was used in Theorem \ref{thm_BacksteppingRelaxed2}.
This is a consequence of the fact that the measure $R_{1}(f,D)$ does not take into account the state constraint for $D$. On the other hand, the measure $R_{2}(f,D)$ takes into account all state constraints (notice that $R_{2}(f,D)$ ``blows up'' when the state tends to a point on the boundary of the state space). In other words, the measure $R_{2}(f,D)$ is a size functional for system \eqref{eq_systemOpenLoop} with (58) and state space $\mathcal{F} \times (\Dmin,\Dmax)$ in the sense of Reference \citenumns{sontag2022remarks}, while the measure $R_{1}(f,D)$ is not a size functional.

\item An interesting special case of this result is the case where the dilution has no upper bound but it must remain positive, namely, $\Dmin = 0, \Dmax = +\infty$. Since $\KarafParA/\Dmax = 1$ and $\KarafParB/\Dmax =1/\Dstar$, 
the transformation \eqref{diffeo2} in this case simplifies to
%\begin{align}
   $ \zeta = \varPhi^{-1}(D) = \ln \left ({D}/{\Dstar} \right )$,
%    \label{diffeo2*}
%\end{align}
and the controller \eqref{eq_IAS1_controller} simplifies to 
\begin{align}
    u(t) =D(t)\left [(\KarafControlGain{1} + \KarafControlGain{2}) (\Dstar - D(t)) + \KarafControlGain{3}\ln \left(\frac{\Dstar}{D(t)}  \left(\frac{y(t)}{\ystar}\right)^{\KarafControlGain{2}} \right) \right], \label{eq_IAS1_controller*}
\end{align}
whereas the state space becomes the ``positive orthant'' $\mathcal{F} \times (0,+\infty)$ and the measure $R_2$ becomes
$R_{2}(f,D)=\max_{a\in [0,A]}\left| \ln (f(a)/\fstar(a)) \right| + \left|\ln \left ({D}/{\Dstar} \right)\right|$.

\item Theorem \ref{thm_ConstrainedDilutionGlobalAsymStability} guarantees Global Asymptotic Stabilization by means of the output feedback law \eqref{eq_IAS1_controller} with exponential convergence rate exactly as Theorem \ref{thm_BacksteppingRelaxed2} does. To see this, notice that estimate \eqref{eq_IAS2_DecayBound} implies that for every $\epsilon >0$ and for every $(f_{0},D_{0}) \in \mathcal{F} \times (\Dmin,\Dmax)$ there exists $T \geq 0$ for which $R_{2}(f[T],D(T)) \leq \exp[- \varphi \left( R_{2}(f_{0},D_{0}) \right) T] \alpha \left( R_{2}(f_{0},D_{0}) \right) \leq \epsilon$. Using the semigroup property and \eqref{eq_IAS2_DecayBound} we get the estimate 
\begin{align}
\begin{array}{l} {R_{2}(f[t],D(t)) \leq \exp[- \varphi \left( \epsilon \right) (t-T)] \alpha \left( R_{2}(f[T],D(T)) \right)} \\ {\leq \exp[- \varphi \left( \epsilon \right) (t-T)] \alpha \left( \alpha \left( R_{2}(f_{0},D_{0}) \right) \right)}
\end{array}
\end{align}
for all $t \geq T$. Combining, we get 
\begin{align}
R_{2}(f[t],D(t)) \leq \exp[- \varphi \left( \epsilon \right) t] \left(1+\frac{\alpha \left( R_{2}(f_{0},D_{0}) \right) }{\epsilon} \right) ^ \frac{\varphi \left( \epsilon \right)}{\varphi \left( R_{2}(f_{0},D_{0}) \right)}  \alpha \left( \alpha \left( R_{2}(f_{0},D_{0}) \right) \right)
\end{align}
for all $t \geq 0$. The above estimate shows a uniform exponential convergence rate $\varphi \left( \epsilon \right) >0$, for every $\epsilon >0$. Consequently, the additional state constraint $D \in (\Dmin,\Dmax)$, does not exclude the possibility of a uniform exponential convergence rate. 

\end{itemize}

The proof of Theorem \ref{thm_ConstrainedDilutionGlobalAsymStability} shows that the controller gains $k_{1}, k_{2}, k_{3}>0$ strongly affect the convergence properties of the solutions of the closed-loop system \eqref{eq_systemOpenLoop} with \eqref{eq_IAS1_controller}. Similarly, the parameters $\Dmin, \Dmax$ strongly affect the convergence properties of the solutions of the closed-loop system \eqref{eq_systemOpenLoop} with \eqref{eq_IAS1_controller}: a tighter state constraint leads to slower convergence and higher overshoots. The behavior of the solutions of the closed-loop \eqref{eq_systemOpenLoop} with \eqref{eq_IAS1_controller} is illustrated by means of representative simulation example in Figure \ref{fig_simulationActDynKarafyllis}. 

Having discussed Theorem \ref{thm_ConstrainedDilutionGlobalAsymStability} in detail, we next provide its proof. 

\begin{proof}[Proof of Theorem \ref{thm_ConstrainedDilutionGlobalAsymStability}]
     Using the transformations \eqref{eq_transformation}, \eqref{eq_adjointEFpi}, \eqref{diffeo1} we get for the closed-loop system \eqref{eq_systemOpenLoop} with \eqref{eq_IAS1_controller}: 
\begin{subequations}
        \label{eq_ClosedLoopKaraf}
            \begin{align}
            \dot \eta(t) &= \frac{\KarafParA \KarafParB(1-\exp[\zeta(t)])}{(\KarafParB+1)(\KarafParB+\exp[\zeta(t)])} := f_1(\zeta(t))\label{eq_etaDynamicsKarafCL}\\
            \dot \zeta(t) &= (\KarafControlGain{1} + \KarafControlGain{2})f_1(\zeta(t)) - \KarafControlGain{3}\left (\zeta(t) - \KarafControlGain{2} \eta(t) -\KarafControlGain{2}v(\psi_t) \right) \label{eq_zetaDynamicsKarafCL}\\
            \psi(t) &= \int_0^A \tilde k(a) \psi(t-a)\, \diff a .\label{eq_IDEdynamicsKarafCL}
        \end{align}
    \end{subequations} 
where $v(\psi_t)$ is defined by \eqref{ggg} and \eqref{eq-vdef}. Again, like in \eqref{eq_functionalPsiV}, we know that the map $t \mapsto v(\psi_t)$ is well-defined and the closed-loop \eqref{eq_ClosedLoopKaraf} locally admits a unique solution. Consider the functional $G(\psi_t)$ defined by \eqref{eq-G} and define the function     
    \begin{equation}
        U_3(\eta,\zeta) = \frac{1}{2}\eta^2 + \frac{\LyapCoeff{1}}{2}(\zeta-\KarafControlGain{2}\eta)^2
        \label{eq_LyapFcnU3}
    \end{equation}
where $\LyapCoeff{1}= \frac{1}{\KarafControlGain{1}\KarafControlGain{2}}$.
%and $\KarafControlGain{2} = \KarafControlGain{2}$. 
Using definition \eqref{eq_LyapFcnU3} and equations \eqref{eq_etaDynamicsKarafCL}, \eqref{eq_zetaDynamicsKarafCL}, 
%and the facts that $\LyapCoeff{1}= \frac{1}{\KarafControlGain{1}\KarafControlGain{2}}$, $\KarafControlGain{2} = \KarafControlGain{2}$, 
we find that
    \begin{align}
        \frac{\diff}{\diff t}U_3(\eta(t),\zeta(t)) &= \frac{1}{\KarafControlGain{2}}\zeta(t) f_1(\zeta(t)) - \LyapCoeff{1}\KarafControlGain{3}(\zeta(t)-\KarafControlGain{2}\eta(t))^2 + \LyapCoeff{1}\KarafControlGain{2}\KarafControlGain{3}(\zeta(t)-\KarafControlGain{2}\eta(t))v(\psi_{t})
        \label{dada1}
    \end{align} 
for all $t \geq 0$ for which the solution of \eqref{eq_ClosedLoopKaraf} exists. Equation \eqref{dada1} in conjunction with the inequality $2 \KarafControlGain{2}(\zeta-\KarafControlGain{2}\eta)v(\psi) \leq (\zeta-\KarafControlGain{2}\eta)^2 + \KarafControlGain{2}^2 v(\psi)^2$, gives for all $t \geq 0$ for which the solution of \eqref{eq_ClosedLoopKaraf} exists:
    \begin{align}
        \frac{\diff}{\diff t}U_3(\eta(t),\zeta(t)) 
        &\leq \frac{1}{\KarafControlGain{2}}\zeta(t) f_1(\zeta(t)) - \frac{1}{2}\LyapCoeff{1}\KarafControlGain{3}(\zeta(t)-\KarafControlGain{2}\eta(t))^2 + \frac{1}{2}\LyapCoeff{1}\KarafControlGain{2}^2\KarafControlGain{3}v(\psi_{t})^2.
        \label{dadaP1}
    \end{align} 
    Using \eqref{eq-vleqG} we obtain from \eqref{dadaP1} the following differential inequality
\begin{align}
        \frac{\diff}{\diff t}U_3(\eta(t),\zeta(t)) 
        &\leq \frac{1}{\KarafControlGain{2}}\zeta(t) f_1(\zeta(t)) - \frac{1}{2}\LyapCoeff{1}\KarafControlGain{3}(\zeta(t)-\KarafControlGain{2}\eta(t))^2 + \frac{\LyapCoeff{1}\KarafControlGain{2}^2\KarafControlGain{3}}{2} \exp[2\sigma A] G(\psi_{t})^2
        \label{dadaP2}
    \end{align} 
for all $t \geq 0$ for which the solution of \eqref{eq_ClosedLoopKaraf} exists. On the other hand, the differential inequality \eqref{eq_BoundPsiLyap} in conjunction with Lemma 2.12 on pages 77-78 in Reference \citenumns{karafyllisJiang2011stability} implies the following estimate
    \begin{equation}
        G(\psi_t) \leq \exp[- \sigma t]  G(\psi_{0})
        \label{dadaP3}
    \end{equation}
for all $t \geq 0$ and consequently, we obtain from \eqref{dadaP2}, \eqref{dadaP3}, the following differential inequality that holds for all $t \geq 0$ for which the solution of \eqref{eq_ClosedLoopKaraf} exists:
    \begin{align}
        \frac{\diff}{\diff t}U_3(\eta(t),\zeta(t)) 
        &\leq \frac{1}{\KarafControlGain{2}}\zeta(t) f_1(\zeta(t)) - \frac{1}{2}\LyapCoeff{1}\KarafControlGain{3}(\zeta(t)-\KarafControlGain{2}\eta(t))^2 + \frac{\LyapCoeff{1}\KarafControlGain{2}^2\KarafControlGain{3}}{2}  \exp[-2\sigma (t-A)] G(\psi_{0})^2.
        \label{dadaP4}
    \end{align}
Notice that by virtue of the fact that $f_1(\zeta)=\frac{\KarafParA \KarafParB(1-\exp[\zeta])}{(\KarafParB+1)(\KarafParB+\exp[\zeta])}$ (which implies that $\zeta f_1(\zeta) \leq 0$ for all $\zeta \in \realNumbers$), we obtain from \eqref{dadaP4} the following estimate for all $t \geq 0$ for which the solution of \eqref{eq_ClosedLoopKaraf} exists:
\begin{align}
        \frac{\diff}{\diff t}U_3(\eta(t),\zeta(t)) 
        &\leq \frac{\LyapCoeff{1}\KarafControlGain{2}^2\KarafControlGain{3}}{2} \exp[-2\sigma (t-A)] G(\psi_{0})^2.
        \label{dadaP5}
    \end{align}
The above estimate in conjunction with Lemma 2.12 on pages 77-78 in Reference \citenumns{karafyllisJiang2011stability} implies that 
\begin{align}
U_3(\eta(t),\zeta(t)) \leq U_3(\eta(0),\zeta(0)) + \frac{\LyapCoeff{1}\KarafControlGain{2}^2\KarafControlGain{3}}{2\sigma}  \exp[2\sigma A] G(\psi_{0})^2
\label{dadaP6}
    \end{align}
for all $t \geq 0$ for which the solution of \eqref{eq_ClosedLoopKaraf} exists. For the positive definite quadratic function $U_3$ defined by \eqref{eq_LyapFcnU3}, there exist constants $\underline{c},\overline{c}>0$ such that for all $(\eta,\zeta)\in \realNumbers^2$ it holds that
    \begin{align}
        \underline c \left( \eta^2 + \zeta ^2\right) \leq U_3(\eta,\zeta) \leq \overline{c}\left( \eta^2 + \zeta ^2\right).
        \label{dada6}
    \end{align}
Inequalities \eqref{dada6}, \eqref{dadaP6} show that the component $(\eta(t),\zeta(t))$ of the solution of \eqref{eq_ClosedLoopKaraf} is bounded for all times $t \geq 0$ for which the solution exists. Thus, the solution of \eqref{eq_ClosedLoopKaraf} exists for all $t \geq 0$ and satisfies 
\begin{align}
\eta^2 (t) + \zeta^2 (t) \leq \frac{\overline{c}}{\underline{c}} \left( \eta^2 (0) + \zeta^2 (0) \right) + \frac{\LyapCoeff{1}\KarafControlGain{2}^2\KarafControlGain{3}}{2\sigma \underline{c}}  \exp[2\sigma A] G(\psi_{0})^2
\label{dadaP7}
    \end{align}
for all $t \geq 0$. The fact that $f_1(\zeta)=\frac{\KarafParA \KarafParB(1-\exp[\zeta])}{(\KarafParB+1)(\KarafParB+\exp[\zeta])}$ and \eqref{eq_LyapFcnU3}, \eqref{dada6} imply the existence of a non-increasing, positive mapping $\varphi(s)>0$ defined for all $s \geq 0$ with the property that for all $s \geq 0$ the following estimate holds for all $(\eta,\zeta)\in \realNumbers^2$ with $\zeta^2 \leq \frac{\overline{c}}{\underline{c}} s^2 + \frac{\LyapCoeff{1}\KarafControlGain{2}^2\KarafControlGain{3}}{2\sigma \underline{c}}  \exp[2\sigma A] s^2$: 
    \begin{align}
        \frac{1}{\KarafControlGain{2}}\zeta f_1(\zeta) - \frac{1}{2}\LyapCoeff{1}\KarafControlGain{3}(\zeta-\KarafControlGain{2}\eta)^2 
        \leq -2 \varphi(s) U_3(\eta,\zeta).
        \label{dadaP8}
    \end{align}
It follows from \eqref{dadaP4} and \eqref{dadaP7}, \eqref{dadaP8} that the following estimate holds for all $t \geq 0$: 
    \begin{align}
        \frac{\diff}{\diff t}U_3(\eta(t),\zeta(t)) 
        &\leq -2 \varphi \left( p_{0} \right)  U_3(\eta(t),\zeta(t)) + \frac{\LyapCoeff{1}\KarafControlGain{2}^2\KarafControlGain{3}}{2}  \exp[-2\sigma (t-A)] G(\psi_{0})^2
        \label{dadaP9}
    \end{align}
where $p_{0}:=|\eta(0)|+|\zeta(0)|+G(\psi_{0})$. Without loss of generality, we may assume that $\varphi(s) \leq \sigma$ for all $s \geq 0$. The differential inequality \eqref{dadaP9} in conjunction with Lemma 2.12 on pages 77-78 in Reference \citenumns{karafyllisJiang2011stability} implies the following estimate for all $t \geq 0$:
    \begin{align}
        \exp[2 \varphi \left( p_{0} \right) t] U_3(\eta(t),\zeta(t)) 
        &\leq  U_3(\eta(0),\zeta(0)) + \frac{\LyapCoeff{1}\KarafControlGain{2}^2\KarafControlGain{3}}{4 \varphi \left( p_{0} \right)}  \exp[2 \sigma A] G(\psi_{0})^2.
        \label{dadaP10}
    \end{align}
 Combining \eqref{dadaP10}, \eqref{dada6}, \eqref{dadaP3}, we get for all $t \geq 0$
    \begin{align}
    \begin{array}{l} {  |\eta (t)| + |\zeta (t)| + \exp[ \sigma A] G(\psi_{t}) } \\
        {\leq  K \left( p_{0} \right) \exp[- \varphi \left( p_{0} \right) t]
        \left( |\eta (0)| + |\zeta  (0)| + \exp[ \sigma A] G(\psi_{0}) \right) }
        \end{array}
        \label{dadaP11}
    \end{align}
where $K(s)$ is the non-decreasing function defined for $s \geq 0$ by the equation
    \begin{align} 
        K(s):=   3\sqrt{\frac{\overline{c}}{\underline{c}}} + \sqrt {\frac{\LyapCoeff{1}\KarafControlGain{2}^2\KarafControlGain{3}}{ \underline{c} \varphi(s)  }}.  
        \label{dadaP12}
    \end{align}
  Using \eqref{eq_KK17_EtaPsiBoundToF_upper}, \eqref{eq_KK17technicalbounds}, \eqref{diffeo2}, \eqref{dadaP12}, the fact that $p_{0}:=|\eta(0)|+|\zeta(0)|+G(\psi_{0})$ and definition \eqref{eq_IAS2_DecayBound22}, we obtain the stability estimate \eqref{eq_IAS2_DecayBound} for an appropriate $\alpha \in \mathcal{K}_\infty$.    The proof is complete. 
\end{proof}

\section{Full-State Stabilization under State Constraints}
\label{sec-constrained-miroslav}

The design that we just presented is mindful of three practical requirements:
\begin{enumerate}
    \item[(i)] dilution must remain above a given positive lower bound $\Dmin >0$;
    \item[(ii)] dilution must obey a given upper bound $\Dmax < +\infty$;
    \item[(iii)] only a bulk concentration $y(t)$ is available for measurement, rather than an age-structured density $f(a,t)$. 
\end{enumerate}
The price we pay for meeting all these practical requirements is a design that is quite complicated and that offers little insight and generalizability to other situations. In the remainder of this section we present an alternative design, which fails to obey two of the above requirements but is clear, elegant, and systematic. 

First, we dispose of the requirement for a strictly positive lower bound $\Dmin >0$ and only require the dilution to remain positive. This is practically reasonable. Dilution/harvesting can be easily completely shut off. Second, we lift the fixed upper bound requirement $\Dmax < +\infty$ and require the dilution to merely be bounded. 

Third, we allow the measurement of the full  state of the system, namely, of the age-structured density $f(a,t)$. This is a stronger requirement as the age distribution of the population cannot be measured. However, the size distribution of the population, which is proportional to the age distribution can be measured. We allow the full-state measurement only for pedagogical reasons, and without loss of generality. The analysis is much cleaner, clearer, and exact when the full $f(a,t)$ is measured, as opposed to when only $y(t)$ is measured. 

So, we start by setting $\Dmin = 0$ and $\Dmax = +\infty$. Then, we recall that
\begin{eqnarray}
    \eta &=& \ln \Pi(f)\\
    \zeta &=& \ln\frac{D}{\Dstar}.
\end{eqnarray}
Next, we introduce a backstepping transformation 
\begin{equation}
    z= \zeta - c_1 \eta, \qquad c_1>0.
\end{equation}
Straightforward calculation then yields
\begin{eqnarray}
\dot\eta &=& \Dstar \left(1- \exp[\zeta]\right)\nonumber\\
&=& \Dstar \left(1- \exp[c_1\eta]\right)+\Dstar\exp[c_1\eta] \left(1- \exp[z]\right)\\
\dot \zeta &=& \frac{u}{D} \\
\dot z &=& \frac{u}{D} - c_1 \Dstar\left(1- \exp[\zeta]\right). 
\end{eqnarray}

Before we design a controller, we introduce two positive definite radially unbounded functions,
\begin{eqnarray}
\omega(z) &=& \exp[z] -1-z\\
\mu(z) &=& \sinh^2\left(\frac{z}{2}\right).
\end{eqnarray}
Next, we introduce two Lyapunov functions,
\begin{eqnarray}
    V_1 &=& \omega(-c_1\eta) \\
    V_2 &=& \omega(\zeta - c_1\eta) = \omega(z)
\end{eqnarray}
and the overall Lyapunov function 
\begin{equation}
    V = \theta V_1 + V_2 = \theta\omega(-c_1\eta) + \omega(\zeta-c_1 \eta), \qquad \theta>0. 
\end{equation}
Noting that $\omega'(z) = \exp[z] -1$, as well as that 
\begin{equation}
\left(\exp[z] -1\right)\left(\exp[-z] -1\right) = - 4 \sinh^2(z/2),
\end{equation}
after a lengthy calculation, with intermediate steps that produce
\begin{eqnarray}
\dot V_1 &= & - 4\Dstar
 c_1 \mu(-c_1\eta) - \Dstar c_1 \left( \exp[z] -1\right)\left( \exp[c_1\eta] -1\right)
\\
\dot V_2 &=& \Dstar\left( \exp[z] -1\right)
\left(\frac{u}{\Dstar D} - c_1 \left(1- \exp[\zeta]\right) \right),
\end{eqnarray}
we note that the feedback law 
\begin{equation}\label{eq-bkst-cont-mk}
u = \Dstar D \left\{
c_1 \left[\theta \left( \exp[c_1\eta] -1\right) + 1- \exp[\zeta]
\right]
+ c_2 \left(\exp[ c_1\eta -\zeta] -1\right)
\right\},
\end{equation}
with $c_2>0$, produces a particularly elegant Lyapunov derivative, 
\begin{equation}
    \dot V = - 4\Dstar\left[
    \theta c_1 \mu(-c_1\eta) + c_2 \mu(\zeta-c_1\eta)
    \right],
\end{equation}
which is negative definite and radially unbounded. Global asymptotic stability of the equilibrium $\eta = \zeta =0$ then follows. With the global exponential stability of the decoupled $\psi$-system, the global asymptotic stability of the $(\eta,\psi,\zeta)$-system follows.

\begin{theorem}
The controller \eqref{eq-bkst-cont-mk}, rewritten in the original $(f,D)$-variables as
\begin{equation}\label{eq-bkst-cont-mk+}
u = \Dstar D \left\{
c_1 \left[\theta \left( \left(\Pi(f)\right)^{c_1} -1\right) + 1- \frac{D}{\Dstar}
\right]
+ c_2 \left(\frac{\Dstar}{D}\left(\Pi(f)\right)^{c_1} -1\right)
\right\}, 
\end{equation}
with gains $c_1,c_2,\theta>0$, guarantees that there exists a class ${\cal KL}$ function $\beta$ such that 
\begin{align}
    %&
    \max_{a\in [0,A]}\left| \ln \frac{f(a,t)}{\fstar(a)} \right| + \left|\ln \frac{D(t)}{\Dstar} \right| \
    %\nonumber\\  &
    \leq 
    \beta\left(\max_{a\in [0,A]}\left| \ln \frac{f_0(a)}{\fstar(a)} \right| + \left|\ln  \frac{D_0}{\Dstar} \right|,t\right), \qquad \forall t\geq 0
\end{align}
for all $f_0\in \mathcal{F}$ and $D_0>0$. 
\end{theorem}
 
This design with $D\in (0, +\infty)$ calls for a reexamination of the more difficult problem where $D\in (\Dmin,\Dmax)$ and $0 \leq \Dmin<\Dstar<\Dmax$. We omit much of the detail and summarize the design and analysis.

Let us denote 
\begin{equation}
    \delta_1 = \frac{\Dstar-\Dmin}{\Dstar}, \qquad n = \frac{\Dstar - \Dmin}{\Dmax-\Dstar}. 
\end{equation}
One gets
\begin{eqnarray}
    \dot\eta &=& \Dstar\delta_1 \frac{1-\exp[\zeta]}{1+n \exp[\zeta]}, \qquad \zeta = \ln\left(\frac{1}{n} \frac{D-\Dmin}{\Dmax-D}\right)\\
    \dot \zeta &=& \frac{\Dmax-D}{(D-\Dmin)(\Dmax-D)}u
\end{eqnarray}
and
\begin{eqnarray}
\dot V_1 &= & - 4\Dstar
 c_1 \frac{\delta_1}{1+n\exp[c_1\eta]}\mu(-c_1\eta) 
 \nonumber\\
 && - \Dstar c_1 \frac{\delta_1(1+n)}{\left(1+n \exp[c_1\eta]\right)^2}\left(\exp[c_1\eta] -1\right)\left( \exp[z] -1\right).
% \\
% \dot V_2 &=& \Dstar\left( \exp[z] -1\right)
% \left(\frac{u}{\Dstar D} - c_1 \left(1- \exp[\zeta]\right) \right),
\end{eqnarray}
With further calculations, one obtains
\begin{equation}\label{eq-Vdot-most-general}
    \dot V = - 4\Dstar\left[
    \theta c_1 \frac{\delta_1}{1+n\exp[c_1\eta]} \mu(-c_1\eta)+ c_2 \mu(\zeta-c_1\eta)
    \right]
\end{equation}
with the controller
\begin{eqnarray}\label{eq-bkst-cont-mk++}
u &=& \frac{\Dstar}{\Dmax-\Dmin}(D-\Dmin)(\Dmax-D) 
\left\{
c_1 \left[\theta \frac{\delta_1(1+n)}{\left(1+n \exp[c_1\eta]\right)^2}\left( \exp[c_1\eta] -1\right) 
+ \frac{1- \exp[\zeta]}{1+n\exp[\zeta]}
\right]
\right. \nonumber\\
&& \left. + c_2 \left(\exp[ c_1\eta -\zeta] -1\right)
\right\}
\nonumber\\
&=&\frac{\Dstar}{\Dmax-\Dmin} (D-\Dmin)(\Dmax-D)\left\{
c_1 \left[\theta \frac{\delta_1(1+n)}{\left(1+n\left(\Pi(f)\right)^{c_1}\right)^2}\left( \left(\Pi(f)\right)^{c_1} -1\right) + 1- \frac{1}{n} \frac{D-\Dmin}{\Dmax-D}
\right]
\right. \nonumber\\
&& \left.+ c_2 \left(
n\frac{\Dmax-D}{D-\Dmin}
\left(\Pi(f)\right)^{c_1} -1\right)
\right\}.  
\end{eqnarray}
This is clearly a considerably more complicated controller than \eqref{eq_IAS1_controller} but the simplicity of the negative definite, radially unbounded \eqref{eq-Vdot-most-general} makes this controller worth consideration.

\begin{theorem}
The controller \eqref{eq-bkst-cont-mk++} with gains $c_1,c_2,\theta>0$, guarantees that there exists a class ${\cal KL}$ function $\beta$ such that 
\begin{align}
    &
    \max_{a\in [0,A]}\left| \ln \frac{f(a,t)}{\fstar(a)} \right| + \left|\ln\left(
    \frac{D(t)-\Dmin}{\Dstar - \Dmin} 
    \frac{\Dmax-\Dstar}{\Dmax-D(t)}
    \right) \right| \
    \nonumber\\  &
    \leq 
    \beta\left(\max_{a\in [0,A]}\left| \ln \frac{f_0(a)}{\fstar(a)} \right| 
    + \left|\ln\left(
    \frac{D_0-\Dmin}{\Dstar - \Dmin} 
    \frac{\Dmax-\Dstar}{\Dmax-D_0}
    \right)\right|,t\right), \qquad \forall t\geq 0
\end{align}
for all $f_0\in \mathcal{F}$ and $D_0\in (\Dmin,\Dmax)$. 
\end{theorem}

\section{Conclusion}

The interest in extending the foundational design and analysis results in Reference \citenumns{KK17Chemostat} goes in many directions. Among them, the most interesting at this stage are in incorporating additional dynamics, which may arise for a variety of reasons (the presence of actuator dynamics, a substrate, input delay, additional species, etc.). In this paper we made what is probably most natural non-trivial step---incorporating actuator dynamics in the form of a single state. Compensating these actuator dynamics was not exceptionally difficult but ensuring that, in addition to the population density remaining positive, the state of the actuator dynamics remains positive as well, or even remains within a given positive interval, is far from elementary and  charts a path towards further generalizations. 

Among the further generalizations are those that bring into the model additional infinite-dimensional states. One such addition is an input delay, which can be compensated using predictor feedback. Other additions of infinite dimensional states are additional species. Even two species can be interconnected in several ways, of which some require no innovations in the control design, while other types of interconnections do. Such generalizations of the chemostat problem into population systems, of which epidemiology is one possible application, open exciting possibilities for future research. 

%\newpage

% \backmatter

%{\section*{Acknowledgments}
%\color{red}This is acknowledgment text. Provide text here. 

\subsection*{Author contributions}

All the authors have contributed to conceiving the problems and writing the article. The simulations were produced by Paul-Erik Haacker.

\subsection*{Financial disclosure}

The work of Miroslav Krstic was funded by the NSF grant ECCS-2151525 and the AFOSR grant FA9550-22-1-0265. The work of Mamadou Diagne was funded by the NSF CAREER Award CMMI-2302030 and  the NSF grant CMMI-2222250.

\subsection*{Conflict of interest}

The authors declare no potential conflict of interests.

%\section*{Supporting information}

%The following supporting information is available as part of the online article:

\appendix
\section{Proof of technical lemmas}
\begin{proof}[Proof of Lemma \ref{lemma_T_1_boundMain}.]We rewrite \eqref{eq_functionalPsi} 
    \begin{align}
        -v_1(\psi_t) = \frac{-1}{\int_0^Ap(a)\fstar(a)(1+\psi(t-a))\diff a} T_3(\psi_t)
    \end{align}
    where
    \begin{align}
        T_3(\psi_t) &= p(A)\fstar(A)(T_4(\psi_t)-\psi(t-A)) - p(0)\fstar(0)(T_4(\psi_t)-\psi(t))\label{eq_def_T3}\\
        &- \int_0^A \tilde p(a) \fstar(a) \diff a T_4(\psi_t) + \int_0^A \tilde p(a) \fstar(a) \psi(t-a)\diff a\\
        T_4(\psi_t) &= \frac{1}{\ystar}\int_0^A p(a)\fstar(a) \psi(t-a)\diff a.
    \end{align}
    By means of the definition of $\ystar$ \eqref{eq_ystar_def} and the nonnegativity of $p$, $\fstar$ on $[0,A]$, it holds that
    \begin{align}
        |T_4(\psi_t)| \leq \max_{a\in[0,A]}|\psi(t-a)|. \label{eq_T4_bound}
    \end{align}
    Taking the absolute value of \eqref{eq_def_T3} and using \eqref{eq_T4_bound} we arrive at
    \begin{align}
        |T_3(\psi_t)| \leq \sqrt{c_1}\ystar\max_{a\in[0,A]}|\psi(t-a)| \label{eq_T3_bound}
    \end{align}
     with
    \begin{align}
       \frac{1}{2} \sqrt{c_1} \ystar= p(A)\fstar(A) + p(0)\fstar(0) + \int_0^A|\tilde p(a)|\fstar(a)\diff a >0.
    \end{align}
    We note, that $c_1$ is independent of the chosen setpoint $\ystar$. With nonnegativity of $p(a)$ and $\fstar(a)$ for $a \in [0,A]$ and $1 + \psi(t-a)>0$ for all $(a,t) \in [0,A] \times \realNumbers_+$, we note that
    \begin{align}
        |v_1(\psi_t)| = \frac{1}{\int_0^Ap(a)\fstar(a)(1+\psi(t-a))\diff a} |T_3(\psi_t)|
    \end{align} and
    \begin{align}
        \int_0^A p(a) \fstar(a) \psi(t-a)\diff a \geq \ystar \min_{a\in[0,A]} \psi(t-a)
    \end{align}
    yields \eqref{eq_T_1_boundMain} using \eqref{eq_T3_bound}.
\end{proof}

% \nocite{*}% Show all bib entries - both cited and uncited; comment this line to view only cited bib entries;

%

% \clearpage

% \section*{Author Biography}
% \color{red}
% \begin{biography}{\includegraphics[width=66pt,height=86pt,draft]{figures/ActuatorDynamics_OpenLoopInstability.eps}}{\textbf{Author Name.} This is sample author biography text this is sample author biography text this is sample author biography text this is sample author biography text this is sample author biography text this is sample author biography text this is sample author biography text this is sample author biography text this is sample author biography text this is sample author biography text this is sample author biography text this is sample author biography text this is sample author biography text this is sample author biography text this is sample author biography text this is sample author biography text this is sample author biography text this is sample author biography text this is sample author biography text this is sample author biography text this is sample author biography text.}
% \end{biography}


\begin{thebibliography}{10}
\providecommand \doibase [0]{http://dx.doi.org/}%

\bibitem{kuritz2017relationship}
Kuritz K, St{\"o}hr D, Pollak N, Allg{\"o}wer F. On the relationship between
  cell cycle analysis with ergodic principles and age-structured cell
  population models. {\it Journal of theoretical biology} 2017\string;
  414\string: 91--102.

\bibitem{billy2012age}
Billy F, Clairambault J, Delaunay F, Feillet CA, Robert N. Age-structured cell
  population model to study the influence of growth factors on cell cycle
  dynamics.. {\it Mathematical Biosciences and Engineering} 2012\string: xx.

\bibitem{rong2007mathematical}
Rong L, Feng Z, Perelson AS. Mathematical analysis of age-structured HIV-1
  dynamics with combination antiretroviral therapy. {\it SIAM Journal on
  Applied Mathematics} 2007\string; 67(3)\string: 731--756.

\bibitem{KK17Chemostat}
Karafyllis I, Krstic M. Stability of integral delay equations and stabilization
  of age-structured models. {\it ESAIM: Control, Optimisation and Calculus of
  Variations} 2017\string; 23(4)\string: 1667--1714.

\bibitem{kurth2021tracking}
Kurth AC, Schmidt K, Sawodny O. Tracking-control for age-structured population
  dynamics with self-competition governed by integro-PDEs. {\it Automatica}
  2021\string; 133\string: 109850.

\bibitem{kurth2023control}
Kurth AC, Sawodny O. Control of age-structured population dynamics with
  intraspecific competition in context of bioreactors. {\it Automatica}
  2023\string; 152\string: 110944.

\bibitem{de2001stabilization2}
De~Leenheer P, Aeyels D. Stabilization of positive linear systems. {\it Systems
  \& control letters} 2001\string; 44(4)\string: 259--271.

\bibitem{de2001stabilization}
De~Leenheer P, Aeyels D. Stabilization results for positive systems with first
  integrals. {\it IFAC Proceedings Volumes} 2001\string; 34(6)\string:
  881--886.

\bibitem{inaba2017age}
Inaba H. {\it Age-structured population dynamics in demography and
  epidemiology}.
\newblock Springer .
\newblock 2017.

\bibitem{martcheva2015introduction}
Martcheva M. {\it An introduction to mathematical epidemiology}. 61.
\newblock Springer .
\newblock 2015.

\bibitem{toth2006limit}
Toth D, Kot M. Limit cycles in a chemostat model for a single species with age
  structure. {\it Mathematical biosciences} 2006\string; 202(1)\string:
  194--217.

\bibitem{m1925applications}
McKendrick AG. Applications of mathematics to medical problems. {\it
  Proceedings of the Edinburgh Mathematical Society} 1925\string; 44\string:
  98--130.

\bibitem{brauer2012mathematical}
Brauer F, Castillo-Chavez C. {\it Mathematical models in population biology and
  epidemiology}. 2.
\newblock Springer .
\newblock 2012.

\bibitem{kermack1927contribution}
Kermack WO, McKendrick AG. A contribution to the mathematical theory of
  epidemics. {\it Proceedings of the Royal Society of London. Series A,
  Containing papers of a mathematical and physical character} 1927\string;
  115(772)\string: 700--721.

\bibitem{boucekkine2011optimalControl}
Boucekkine R, Hritonenko N, Yatsenko Y. {\it Optimal control of age-structured
  populations in economy, demography, and the environment}.
\newblock Routledge .
\newblock 2011.

\bibitem{charlesworth1994ageStructured}
Charlesworth B. {\it Evolution in age-structured populations}. 2.
\newblock Cambridge University Press Cambridge .
\newblock 1994.

\bibitem{rundnicki1994asymptoticSimilarity}
Rundnicki R, Mackey MC. Asymptotic similarity and Malthusian growth in
  autonomous and nonautonomous populations. {\it Journal of Mathematical
  Analysis and Applications} 1994\string; 187(2)\string: 548--566.

\bibitem{smith1995theoryOfChemostat}
Smith HL, Waltman P. {\it The theory of the chemostat: dynamics of microbial
  competition}. 13.
\newblock Cambridge university press .
\newblock 1995.

\bibitem{inaba1988semigroupErgodic}
Inaba H. A semigroup approach to the strong ergodic theorem of the multistate
  stable population process. {\it Mathematical Population Studies} 1988\string;
  1(1)\string: 49--77.

\bibitem{inaba1988asymptoticProperties}
Inaba H. Asymptotic properties of the inhomogeneous Lotka-Von Foerster system.
  {\it Mathematical Population Studies} 1988\string; 1(3)\string: 247--264.

\bibitem{feichtinger2003optimality}
Feichtinger G, Tragler G, Veliov VM. Optimality conditions for age-structured
  control systems. {\it Journal of Mathematical Analysis and Applications}
  2003\string; 288(1)\string: 47--68.

\bibitem{sun2014optimal}
Sun B. Optimal control of age-structured population dynamics for spread of
  universally fatal diseases II. {\it Applicable Analysis} 2014\string;
  93(8)\string: 1730--1744.

\bibitem{schmidt2018trajectory}
Schmidt K, Karafyllis I, Krstic M. Yield trajectory tracking for hyperbolic
  age-structured population systems. {\it Automatica} 2018\string; 90\string:
  138--146.

\bibitem{karafyllis2008controlODE}
Karafyllis I, Kravaris C, Syrou L, Lyberatos G. A vector Lyapunov function
  characterization of input-to-state stability with application to robust
  global stabilization of the chemostat. {\it European Journal of Control}
  2008\string; 14(1)\string: 47--61.

\bibitem{karafyllis2009relaxedODEcontrol}
Karafyllis I, Kravaris C, Kalogerakis N. Relaxed Lyapunov criteria for robust
  global stabilisation of non-linear systems. {\it International Journal of
  Control} 2009\string; 82(11)\string: 2077--2094.

\bibitem{gouze2006robust}
Gouze JL, Robledo G. Robust control for an uncertain chemostat model. {\it
  International Journal of Robust and Nonlinear Control} 2006\string;
  16(3)\string: 133--155.

\bibitem{de2003feedback}
De~Leenheer P, Smith H. Feedback control for chemostat models. {\it Journal of
  Mathematical Biology} 2003\string; 46(1)\string: 48--70.

\bibitem{dimitrova2012nonlinear}
Dimitrova N, Krastanov M. Nonlinear adaptive stabilizing control of an
  anaerobic digestion model with unknown kinetics. {\it International Journal
  of Robust and Nonlinear Control} 2012\string; 22(15)\string: 1743--1752.

\bibitem{sontag2022remarks}
Sontag ED. Remarks on input to state stability of perturbed gradient flows,
  motivated by model-free feedback control learning. {\it Systems \& Control
  Letters} 2022\string; 161\string: 105138.

\bibitem{karafyllisJiang2011stability}
Karafyllis I, Jiang ZP. {\it Stability and stabilization of nonlinear systems}.
\newblock Springer Science \& Business Media .
\newblock 2011.

\bibitem{toth2008strong}
Toth D. Strong resonance and chaos in a single-species chemostat model with
  periodic pulsing of resource. {\it Chaos, Solitons \& Fractals} 2008\string;
  38(1)\string: 55--69.

\bibitem{toth2008bifurcation}
Toth D. Bifurcation structure of a chemostat model for an age-structured
  predator and its prey. {\it Journal of Biological Dynamics} 2008\string;
  2(4)\string: 428--448.

\bibitem{kurth2021optimalControl}
Kurth AC, Schmidt K, Sawodny O. Inversion-based and optimal feedforward control
  for population dynamics with input constraints and self-competition in
  chemostat reactor applications. {\it Journal of Dynamic Systems, Measurement,
  and Control} 2021\string; 143(5).

\bibitem{schmidt2017Simulation}
Schmidt K, Sawodny O. Efficient simulation of semilinear populations models for
  age-structured bio reactors. In: IEEE. ; 2017\string: 1716--1721.

\bibitem{10130642}
Krstic M. Inverse Optimal Safety Filters. {\it IEEE Transactions on Automatic
  Control} 2023\string: 1-16.
\newblock \href {\doibase 10.1109/TAC.2023.3278788} {doi:
  10.1109/TAC.2023.3278788}

\end{thebibliography}
\end{document}